\documentclass [twoside,reqno, 12pt] {amsart}

\usepackage{fullpage}

\usepackage{amsfonts,mathtools}
\usepackage{amssymb, tikz}
\usepackage{titlesec}
\usepackage{color}
\usepackage{hyperref}

\usepackage{comment}
\usepackage{caption}
\usepackage{listings}
\usepackage{xcolor}
\usepackage{rotating}
\usepackage{subcaption}
\usepackage{biblatex}
\allowdisplaybreaks

\usepackage{color}    

\makeatletter
\providecommand{\@secnumpunct}{.\enspace}
\makeatother

\newcommand{\R}{\mathbb{R}}
\newcommand{\C}{\mathbb{C}}
\newcommand{\CP}{\mathbb{CP}}
\newcommand{\RP}{\mathbb{RP}}
\renewcommand{\P}{\mathbb{P}}
\newcommand{\A}{\mathbb{A}}
\newcommand{\N}{\mathbb{N}}

\DeclareMathOperator{\disc}{disc}

\DeclareMathOperator{\im}{Im}
\newcommand{\p}{\partial}
\newcommand{\mC}{\mathbb{C}} 
\DeclareMathOperator{\supp}{supp}
\DeclareMathOperator{\Spec}{Spec}
\newcommand{\vev}[1]{\left\langle#1\right\rangle}
\newcommand{\K}{\mathbb{K}}
\newcommand{\spoly}[1]{\sigma_{#1}}
\newcommand{\possymsym}[1]{\operatorname{Sym}_+^2(\operatorname{Sym}^2(#1))}

\newcommand{\norm}[1]{\ensuremath{\left\lVert  #1\right\rVert}} 
\newcommand{\abs}[1]{\ensuremath{\left\vert  #1\right\vert}}

\definecolor{dkgreen}{rgb}{0,0.6,0}
\definecolor{dred}{rgb}{0.545,0,0}
\definecolor{dblue}{rgb}{0,0,0.545}
\definecolor{lgrey}{rgb}{0.9,0.9,0.9}
\definecolor{gray}{rgb}{0.4,0.4,0.4}
\definecolor{darkblue}{rgb}{0.0,0.0,0.6}
\definecolor{turquoise}{rgb}{0.2500,0.8750,0.8125}
\definecolor{indigo}{rgb}{0.2930, 0,0.5078}
\definecolor{mag}{rgb}{1, 0,1}
\definecolor{corn}{rgb}{0.3906,0.5820,0.9258}
\definecolor{mvr}{rgb}{0.7773,0.0820,0.5195}
\definecolor{dod}{rgb}{0.11719,0.5625,1}
\lstdefinelanguage{MatlabCostum}{
      backgroundcolor=\color{white},  
      basicstyle=\footnotesize \ttfamily \color{black} \bfseries,   
      breakatwhitespace=false,       
      breaklines=true,               
      captionpos=b,                   
      commentstyle=\color{dkgreen},   
      emph={repmat, ones},			
      keywordstyle=\color{blue},           
      escapeinside={\%*}{*)},                  
      frame=single,                  
      language=Matlab,                
      identifierstyle=\color{black},
      stringstyle=\color{blue},      
      numbers=left,                 
      numbersep=5pt,                  
      numberstyle=\tiny\color{black}, 
      rulecolor=\color{black},        
      showspaces=false,               
      showstringspaces=false,        
      showtabs=false,                
      stepnumber=1,                   
      tabsize=4,
      title=\lstname
}
\newtheorem{theorem}{Theorem}[section]
\newtheorem{corollary}[theorem]{Corollary}
\newtheorem{definition}[theorem]{Definition}
\newtheorem{lemma}[theorem]{Lemma}
\newtheorem{proposition}[theorem]{Proposition}
\newtheorem{remark}[theorem]{Remark}

\newtheorem{Example}[theorem]{Example}

\newcommand{\joonas}[1]{\footnote{Joonas: #1}}

\allowdisplaybreaks

\title{Boundary determination from the anisotropic elastic Dirichlet-to-Neumann map}

\author[Ilmavirta]{Joonas Ilmavirta}
\address{J. Ilmavirta, Department of Mathematics and Statistics, University of Jyv\"askyl\"a, Jyv\"askyl\"a, 40014,  Finland} \email{joonas.ilmavirta@jyu.fi}

\author[Salo]{Mikko Salo}
\address{M. Salo, Department of Mathematics and Statistics, University of Jyv\"askyl\"a, Jyv\"askyl\"a, 40014,  Finland}
\email{mikko.j.salo@jyu.fi}

\author[Schl\"uter]{Hjørdis Schl\"uter}
\address{H. Schl\"uter, Department of Mathematics and Statistics, University of Helsinki, 00014, Finland} \email{hjordis.schluter@helsinki.fi }

\author[Windisch]{Daniel Windisch}
\address{D. Windisch, Department of Computer Science, KU Leuven, 3001 Leuven, Belgium.}
\email{daniel.windisch.math@gmail.com}

\date{\today}

\begin{document}

\begin{abstract}
We show that the Dirichlet-to-Neumann map of the elastic wave equation determines the density and the anisotropic stiffness tensor on the boundary for generic density and stiffness tensor fields. The proof builds on an interplay between microlocal analysis and algebraic geometry. Microlocal analysis yields a construction of special solutions, valid for many real principal type systems, that concentrate near a boundary point. Generic conditions ensure that this construction holds and allow for extraction of information about stiffness tensor and density at the boundary. We also describe the generic sets through algebraic conditions.
\end{abstract}

\maketitle

\section{Introduction} \label{sec_intro}
We study an anisotropic elastic medium and address the inverse problem of recovering its material parameters, stiffness and density of mass, from boundary measurements. These measurements correspond to imposing a displacement that produces a wave and measuring the corresponding stress at the boundary of the medium. This problem has applications in seismic imaging where one aims at creating images of the subsurface of the earth or other planets.

Formally, we address the inverse problem of recovering an anisotropic stiffness tensor~$c$ and a density of mass~$\rho$ from the Dirichlet-to-Neumann map~$\Lambda_{c,\rho}$ at the boundary of an $n$-dimensional domain $\Omega$ with $n\geq 2$. In particular, we show that $\Lambda_{c,\rho}$ determines a generic pair $(c,\rho)$ uniquely at the boundary. The generic conditions are discussed further below and described in Lemma~\ref{lem:genPair}.\par
Boundary determination for this setting has been obtained for isotropic stiffness tensors in~\cite{Rachele001}. The boundary determination step of results has been generalized to transversely isotropic and orthorhombic stiffness tensors \cite{HoopNakamuraZhai19}. In the stationary case boundary determination for the Lamé system has been obtained in~\cite{Nakamura1993}. Boundary determination results for wave equations can be found in~\cite{SU90,SU05,Mo14,SY18}, and for second order elliptic equations boundary determination goes back to~\cite{KV84,SU88}, with~\cite{LU89,DKSU09} treating the anisotropic geometric case.\par

 For anisotropic stiffness tensors this inverse problem is still poorly understood and our result is a first step towards proving the following interior uniqueness result conjectured in \cite[Conjecture 2.4]{Ilmavirta2025}: $\Lambda_{c,\rho}$ determines generic~$c$ and~$\rho$ uniquely in~$\Omega$.

An anisotropic stiffness tensor is an element of $\possymsym{\R^n}$; these symmetries  and the positivity are described in~\eqref{eq:sym} and~\eqref{eq:pos}, respectively.
The rigorous version of our main result is given as follows for smooth material parameters $\rho,\tilde \rho, c$ and $\tilde c$:


\begin{theorem}\label{thm:main}
Let $\Omega \subset \R^n$ with $n\geq 2$ be a bounded and connected domain with a smooth boundary $\partial \Omega$.
Fix a boundary point $p\in\partial\Omega$.
There is such an open and dense subset $G_p\subset (\R_+\times\possymsym{\R^n})^2$ that for all $(\rho,c,\tilde\rho,\tilde c)\in C^\infty(\bar\Omega;\R_+\times\possymsym{\R^n})^2$ with $(\rho(p),c(p),\tilde\rho(p),\tilde c(p))\in G_p$ the following holds:
if $\Lambda_{\rho,c}=\Lambda_{\tilde \rho,\tilde c}$ with the DN maps defined in~\eqref{eq:DN}, 
then
\begin{equation}
c(p)=\tilde{c}(p) \quad \text{and} \quad \rho(p)=\tilde{\rho}(p)
.
\end{equation}
\end{theorem}


\begin{corollary}
\label{cor:fields}
Let $\Omega \subset \R^n$ with $n\geq 2$ be a bounded and connected domain with a smooth boundary $\partial \Omega$.
There is a residual set $G\subset C^\infty(\partial\Omega;\R_+\times\possymsym{\R^n})^2$ so that if $(\rho,c,\tilde\rho,\tilde c)$ are smooth density and stiffness tensor fields whose boundary restriction is in~$G$ and $\Lambda_{\rho,c}=\Lambda_{\tilde \rho,\tilde c}$, then
\begin{equation}
c=\tilde{c} \quad \text{and} \quad \rho=\tilde{\rho}
\end{equation}
on~$\partial \Omega$.
\end{corollary}

For simplicity we have only given results on determining the boundary values of the unknown coefficients. The question of determining higher order normal derivatives at the boundary will be left to future work.

\subsection{Anisotropic elasticity}
The elastic medium is modeled by $\Omega \subset \R^n$ with $n\geq 2$. The stiffness of the material is modeled as a rank 4 tensor field $c(x)\in \R^{n\times n \times n \times n}$ and the density of mass is modeled as a scalar-valued function $\rho \in \R$.
In the following we make the natural assumption that~$\rho$ is positive:
\begin{equation}\label{eq:rhopos}
    \rho > 0.
\end{equation}
The stiffness tensor links the stress and strain tensors of the material. By the symmetries of stress and strain and energy considerations, we assume that~$c$ has the following minor and major symmetries 
    \begin{equation}\label{eq:sym}
        c_{ijk\ell}=c_{jik\ell}=c_{k\ell ij}.
    \end{equation}
Additionally, we assume that~$c$ is positive definite:\footnote{Here and henceforth we follow the Einstein sum convention despite all our indices being down. That is, $c_{ijkl}A_{ij}A_{kl}\coloneqq\sum_{i,j,k,l=1}^n c_{ijkl}A_{ij}A_{kl}$.}
    \begin{equation}\label{eq:pos}
        c_{ijk\ell}A_{ij}A_{k\ell} \geq \delta A_{ij}^2
    \end{equation}
for some $\delta>0$ and all symmetric matrices~$A$.
The waves are modelled by the following linear elastic wave equation in $\Omega \times (0,T)$:
\begin{equation}\label{Pop}
    (P_{\rho,c}u(x,t))_i := \rho(x) \partial_t^2 u_i(x,t) - \partial_{x_j}(c_{ijk\ell}(x) \partial_{x_k}u_{\ell}(x,t))=0, \qquad i=1,...,n,
\end{equation}
where we denote by $P_{\rho,c}$ the elastic wave operator that acts on a displacement $u\in \R^n$. Imposing a displacement at the boundary of the space time cylinder $\partial \Omega \times [0,T]$ gives rise to the following boundary value problem for $u$:
\begin{equation}\label{eq:EWE}
    \begin{cases}
    P_{\rho,c}u = 0 & \text{in } \Omega \times (0,T)\\
    u = f & \text{on } \partial \Omega \times (0,T)\\
    u(x,0)=\partial_t u(x,0) = 0 & \text{in }\Omega,
    \end{cases}
\end{equation}
with $f(x,0)=0$ and $\partial_t f(x,0) = 0$ for $x \in \partial \Omega$. Assuming that $c_{ijk\ell} \in L^{\infty}(\Omega), \rho \in L^{\infty}(\Omega)$ and $f \in L^2(0,T;H^{\frac{1}{2}}(\partial \Omega)) \cap H^1(0,T;L^2(\partial \Omega))$ the initial value boundary value problem \eqref{eq:EWE} has the unique solution  $u \in L^2(0,T;H^1(\Omega))\cap H^1(0,T;L^2(\Omega))$ and additionally, $u \in C(0,T;H^1(\Omega))\cap C^1(0,T;L^2(\Omega))$ \cite[Thm. 2.4.5]{stolk00}. The boundary measurements that we consider are encoded by the Dirichlet-to-Neumann (DN) map $\Lambda$ that maps the Dirichlet boundary function $f$ to the corresponding Neumann boundary function and which is defined as follows:
\begin{equation}\label{eq:DN}
        (\Lambda_{\rho, c} f)_i=\nu_j \, c_{ijk\ell} \, \partial_{x_{k}} u_{\ell}  \Bigg\vert_{(0,T) \times \partial \Omega}, \qquad i=1,...,n
\end{equation}
where $\nu$ denotes the unit outward normal to~$\partial \Omega$.

If we expand the elastic wave operator~$P_{\rho,c}$ from above, we obtain:
\begin{equation}
    (P_{\rho,c})_{i\ell} = \rho\delta_{i\ell} \partial_t^2 - c_{ijk\ell} \partial_j \partial_k - (\partial_j c_{ijk\ell}) \partial_k.
\end{equation}
Fix a point $(x,t)\in\overline{\Omega}\times(0,T)$. We denote a cotangent vector at $(x,t)$ by
\begin{equation}
(\xi,\xi_0)\in T^*_{(x,t)}(\R^n\times\R)\cong\R^n\times\R,
\end{equation}
where $\xi=(\xi_1,\dots,\xi_n)$ is dual to $x$ and $\xi_0$ is dual to $t$.
This yields the full symbol:
\begin{equation}
    (p_{\rho,c})_{i\ell}(x,t,\xi, \xi_0) = \underbrace{- \rho \delta_{i\ell} \xi_0^2 + c_{ijk\ell} \xi_j \xi_k}_{=(p_2)_{i\ell}} \underbrace{- i (\partial_j c_{ijk\ell}) \xi_k}_{=(p_1)_{i\ell}}
\end{equation}
Where we denote the principal symbol by $p_2$ and the lower order term by $p_1$. We define the Christoffel matrix:
\begin{equation}
    (\Gamma_{c})_{i\ell}(x,\xi):= c_{ijk\ell}(x)\xi_j \xi_k
\end{equation}
Positive definiteness of $c$ implies positive definiteness of $\Gamma_c$: Inserting the symmetric matrix $A_{ij}=\frac12(v_i\xi_j+v_j\xi_i)$ in~\eqref{eq:pos} and using~\eqref{eq:sym} gives
\begin{equation}\label{eq:Gamma-elliptic}
v^T\Gamma_c(x,\xi)v=c_{ijk\ell}(x)A_{ij}A_{k\ell}\geq\delta A_{ij}A_{ij}\geq\tfrac{\delta}{2}\abs{v}^2\abs{\xi}^2
\end{equation}
for $v,\xi\in\R^n$, such that $\Gamma_c(x,\xi)$ is positive definite for $\xi\neq0$.
With this definition the principal symbol is on the form:
\begin{equation}
    p_2(x,\xi)=\Gamma_c(x,\xi) - \rho(x) \xi_0^2 I.
\end{equation}
A solution to the elastic wave equation can only have a singularity at $(x,t,\xi,\xi_0)$ if 
    \begin{equation}
        [\Gamma_c(x,\xi)-\rho(x)\xi_0^2 I]\omega = 0,
    \end{equation}
where $\omega$ is an eigenvector of the Christoffel matrix corresponding to an eigenvalue $\lambda(\xi)=\rho \xi_0^2$. Physically, $\xi$ has the interpretation of momentum and $\omega$ has the interpretation of polarization. The so-called \emph{slowness surface} describes the vanishing set of the principal symbol in $T^*_{(x,t)}(\R^n\times\R)$ and thus the kinematics of singularities:
\begin{equation}
        \Sigma_{\rho,c} = \{(\xi,\xi_0) \in T^*_{(x,t)}(\R^n\times\R); \det(\Gamma_c(x,\xi)-\rho(x)\xi_0^2 I)=0\}.
    \end{equation}
This is the vanishing set of the slowness polynomial
\begin{equation}
    \spoly{\rho,c}(x,\xi) := \det \left( \Gamma_c(x,\xi) - \rho \xi_0^2 I \right).
\end{equation}
The slowness polynomial can be written in terms of the eigenvalues of the Christoffel matrix $\lambda_i(x,\xi)$ for $i=1,\dots,n$:
\begin{equation}
    \spoly{\rho,c}(x,\xi)=(\lambda_1(x,\xi)-\rho(x) \xi_0^2)(\lambda_2(x,\xi)-\rho(x) \xi_0^2)\cdots (\lambda_n(x,\xi)-\rho(x) \xi_0^2)
\end{equation}
such that $(\xi,\xi_0) \in \Sigma_{\rho,c}$ if and only if $\lambda_i(x,\xi)=\rho(x) \xi_0^2$ for some $i=1,\dots,n$. Each eigenvalue corresponds to a branch in the slowness surface. In the following we consider polynomials~$\spoly{\rho,c}$ that are square free in the sense that~$\spoly{\rho,c}$ cannot be factored into polynomials that have power larger than 1 and this is equivalent to~$\Gamma_c$ having distinct eigenvalues for some~$\xi$.

\subsection{Idea behind the proof of Theorem~\ref{thm:main}}\label{sec:proofidea}

The proof of Theorem~\ref{thm:main} builds on an interplay between microlocal analysis and algebraic geometry. The algebraic approach builds on ideas in \cite{dehoop2023reconstruction} that address the recovery of an anisotropic stiffness tensor from the slowness surface, and more broadly extend the recent use of algebraic geometry for inverse problems~\cite{maxwell-ag, elastic-regularity}. The first part of the proof uses microlocal analysis and constructs special solutions in the elliptic region concentrating around a boundary point $\hat{x}\in \partial \Omega$ for real principal type operators with matrix valued principal symbol. This generalizes the corresponding result for scalar real principal type operators in \cite[Thm 5.1]{Oksanen2024}.

We present the boundary determination results for general real principal type operators on compact $n$-manifolds with boundary, for which the elastic wave operator in an $n$-dimensional Euclidean domain is just a special case. In the general case, we present what boundary determination results are possible from the Cauchy data set. However, in order to obtain boundary determination results from the Dirichlet-to-Neumann map we specialize to the elastic wave operator in $\Omega$ and the following discussion is centered around this setting.

In order for the construction of special solutions to apply for the elastic operator, we need a covector $\xi'=(\xi_0,\dots,\xi_{n-1})$ such that the roots in~$\xi_n$ of the slowness polynomial~$\spoly{\rho,c}$ are simple and non-real. For positive definite stiffness tensors such covectors exist as soon as~$\spoly{\rho,c}$ is square free (see the proof of Lemma~\ref{lem:alg-alessandrini}). Using the Alessandrini identity it follows from this construction that, for a pair $((\rho,c),(\tilde \rho, \tilde c))$ that satisfies $\Lambda_{\rho,c}=\Lambda_{\tilde \rho, \tilde c}$, the following holds around a boundary point $\hat{x}\in \partial \Omega$:
\begin{equation}\label{eq:Fintro}
    F(\xi,\tilde\xi,\omega,\tilde\omega)
:=
[
(c(\hat{x})-\tilde c(\hat{x}))_{ijkl}
\xi_k\tilde\xi_j
-
(\rho(\hat{x})-\tilde\rho(\hat{x}))
\delta_{il}
\xi_0\tilde\xi_0
]
\omega_l
\tilde\omega_i=0,
\end{equation}
where $\xi,\tilde\xi,\omega,\tilde\omega$ are such that the following eigenvector equations hold:
\begin{align}
    f_i(\xi,\tilde\xi,\omega,\tilde\omega)
&:=
(
c_{ijkl}(\hat{x})\xi_j\xi_k
-
\rho(\hat{x})\delta_{il}\xi_0^2
)
\omega_l=0
,\\
\tilde f_i(\xi,\tilde\xi,\omega,\tilde\omega)
&:=
(
\tilde c_{ijkl}(\hat{x})\tilde\xi_j\tilde\xi_k
-\tilde\rho(\hat{x})\delta_{il}\tilde\xi_0^2
)
\tilde\omega_l=0,
\end{align}
with $i=1,\dots, n$ and where $\xi$ and $\tilde \xi$ satisfy $\xi_i=\tilde \xi_i$ for $i=1,\dots, n-1$. The question is, whether $F=0$ implies $\rho(\hat{x})=\tilde\rho(\hat{x})$ and $c(\hat{x})=\tilde c(\hat{x})$? It seems that this result may not be easy to obtain along the microlocal analysis route. For that purpose, we use algebraic geometry to obtain the desired conclusion.

Motivated by the construction of special solutions, we denote
\begin{equation}\label{eq:Xidef}
\Xi
=
\{
(\xi,\tilde\xi)
\in(\CP^n)^2
;
\xi_i=\tilde\xi_i
\text{ for all }
i=0,\dots,n-1
\}
\end{equation}
and define the following varieties in the variables $\xi, \tilde \xi, \omega, \tilde \omega$:
\begin{equation}\label{eq:Wmeq}
\begin{split}
W_m
&=
\{
(\xi,\tilde\xi,\omega,\tilde\omega)\in \Xi \times (\C \P^n)^2
; \quad
f_i=0
\text{ and }
\tilde f_i=0
\text{ for all }
i=1,\dots,n
\}
\end{split}
\end{equation}
and
\begin{equation}\label{eq:Umeq}
U_m
=
\{
(\xi,\tilde\xi,\omega,\tilde\omega) \in \Xi \times (\C \P^n)^2
;
F
=
0
\}
.
\end{equation}
Algebraic geometry then allows us to conclude that if for $m=(\xi,\tilde \xi, \omega,\tilde \omega) \in \Xi \times (\C \times \C^N)^2$ no irreducible component of~$W_m$ is contained in~$U_m$, it follows that $F=0$ with $f_i=\tilde f_i=0$ implies $\rho(\hat{x})=\tilde\rho(\hat{x})$ and $c(\hat{x})=\tilde c(\hat{x})$. This result holds for generic pairs $(\rho,c)$ and $(\tilde \rho, \tilde c)$ such that their slowness surfaces are square free and no irreducible component of $W_m$ is contained in $U_m$. 
The failure set of squarefreeness we make explicit, while the failure set corresponding to no irreducible component of~$W_m$ being contained in~$U_m$ we do not determine explicitly; this is subject to future work.

\subsection{Structure of the paper} 
This paper is structured as follows: In section \ref{sec:Main} we formalize the results that we need to prove the main result,  Theorem~\ref{thm:main}. Section \ref{sec:specialSols} outlines the construction of special solutions for general real principal type operators with matrix-valued principal symbol on compact $n$-manifolds with boundary. The elastic wave operator in an $n$-dimensional Euclidean domain is a special case for which this construction allows to prove one aspect of the main result. In section~\ref{sec:bdCauchy} we show boundary determination results from Cauchy data for general operators and in section~\ref{sec:bdElastic} we specialize these to boundary determination results from the Dirichlet-to-Neumann map for the elastic wave operator. In section \ref{sec:algebra} we present the algebraic results that are needed to prove the remaining aspects of the main result.  
In section~\ref{sec:fields} we promote the pointwise result to density and stiffness tensor fields.

\subsection*{Declaration of tools}

During the preparation of this manuscript, Anthropic’s
Claude Opus 5.5 was used for minor assistance in the proof of Lemma~\ref{lem:simple-roots-sqf} and Lemma~\ref{lem:squarefreeness}, improving the result of Lemma~\ref{lemma_z_omega} and making the proof of Lemma~\ref{lem:alg-alessandrini} explicit and making the proof of Theorem~\ref{thm:main} explicit once Theorem~\ref{thm:generic-uniqueness} is established. Additionally, Maple 2025 was used to generate Figure~\ref{fig:genpair1} and Figure~\ref{fig:genpair2} and Anthropic’s
Claude Opus 5.5 assisted in writing the corresponding Maple code. LLM tools were also used to help with proofreading the manuscript. The authors take responsibility for all mathematical statements and proofs presented in this paper.

\subsection*{Acknowledgments}
We thank Maria Montanucci and Jonathan Tilling Niemann at the Technical University of Denmark for useful discussions on algebraic aspects of this paper. JI, MS and HS were supported by the Research Council of Finland (Flagship of Advanced Mathematics for Sensing Imaging and Modelling, grant 359208, and Centre of Excellence of Inverse Modelling and Imaging, grants 353091 and 353092). JI was supported by Research Council of Finland grants 351665, 351656, 358047, 360434 and a Väisälä project grant by the Finnish Academy of Science and Letters.
MS was supported by the European Union (ERC AdG 101266232, NETCaP). DW was partially supported by the FWO grants G0F5921N (Odysseus) and G023721N, and by the KU Leuven grant iBOF/23/064. Views and opinions expressed are however those of the author(s) only and do not necessarily reflect those of the European Union or the European Research Council. Neither the European Union nor the granting authority can be held responsible for them.

\section{Idea of proof}
\label{sec:Main}

All material parameters are assumed to be~$C^\infty$ throughout the paper.

%

The explicit generic conditions imposed on a pair $((\rho,c),(\tilde \rho, \tilde c))$ in Theorem~\ref{thm:main} are given as follows:

\begin{lemma}\label{lem:genPair}
    The following conditions on a generic pair $m = ((\rho,c),(\tilde \rho, \tilde c))$ as in Theorem~\ref{thm:main} can be made explicit:
    \begin{enumerate}
        \item[(1)] $\spoly{\rho,c}(\xi_0,\xi_1,\dots,\xi_{n-1},z)$ and $\spoly{\tilde \rho,\tilde c}(\xi_0,\xi_1,\dots,\xi_{n-1},\tilde z)$ are square-free polynomials in\par $\C[\xi_0,\xi_1,\dots,\xi_{n-1},z]$ and $\C[\xi_0,\xi_1,\dots,\xi_{n-1},\tilde z]$ respectively
        \item[(2)] for $m \in (\C \times \C^N)^2$ no irreducible component of $W_m$ is contained in $U_m$
    \end{enumerate}
    and these conditions hold on an open and dense set of pairs that satisfy \eqref{eq:rhopos}--\eqref{eq:pos}.
\end{lemma}

\begin{remark}\label{rem:failureset}
The failure set of squarefreeness is explicit (see Lemma~\ref{lem:generic-simplicity}). For the condition that no irreducible component of $W_m$ lies in $U_m$, we only know that the failure set is contained in a proper algebraic subset of $(\C \times \C^N)^2$ (see Lemma~\ref{lem:generic-noncontainment}), which we do not determine explicitly.
\end{remark}

\begin{figure}[ht]
    \centering
    \begin{minipage}[t]{0.3\textwidth}
        \centering
        \includegraphics[width=0.95\textwidth]{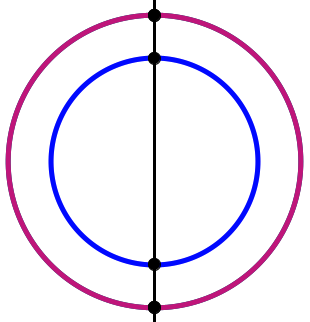}
        \subcaption{$s=0$}
    \end{minipage}%
    \begin{minipage}[t]{0.3\textwidth}
        \centering
        \includegraphics[width=0.95\textwidth]{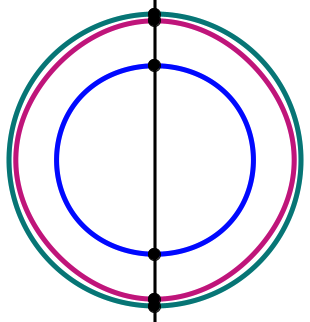}
        \subcaption{$s=0.1$}
    \end{minipage}%
    \begin{minipage}[t]{0.3\textwidth}
        \centering
        \includegraphics[width=0.95\textwidth]{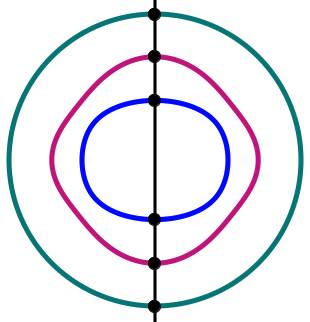}
        \subcaption{$s=1$}
    \end{minipage}
    \caption{The slowness surfaces for the family $c^{(s)}$ when $n=3$ in the $(\xi_2,\xi_3)$-plane of Example~\ref{ex:genpair2}(1). Left: the isotropic tensor ($\lambda=0$, $\mu=1$), for which two branches coincide. Right: the tensor $c^*$. The ray in direction $d=(0,0,\pm 1)$ meets the coinciding branches in a single point on the left and all branches in distinct points on the right. Perturbing in $s$ separates the branches: $\Gamma_{c^*}(d)$ has three distinct eigenvalues for $d=(0,0,1)$. }
    \label{fig:genpair1}
\end{figure}

\begin{figure}[ht]
    \centering
    \begin{minipage}[t]{0.5\textwidth}
        \centering
        \includegraphics[width=0.8\textwidth]{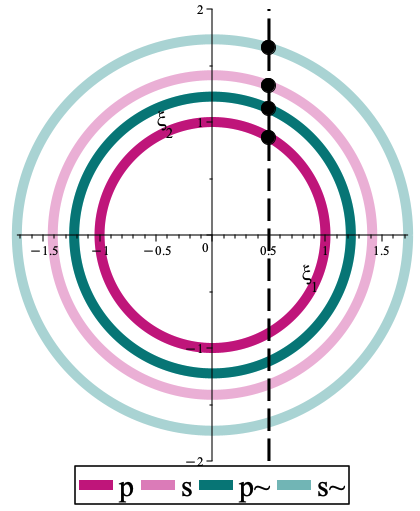}
        \subcaption{$\Sigma_{\rho,c}\cup \Sigma_{\tilde \rho,\tilde c}$}
    \end{minipage}%
    \begin{minipage}[t]{0.5\textwidth}
        \centering
        \includegraphics[width=\textwidth]{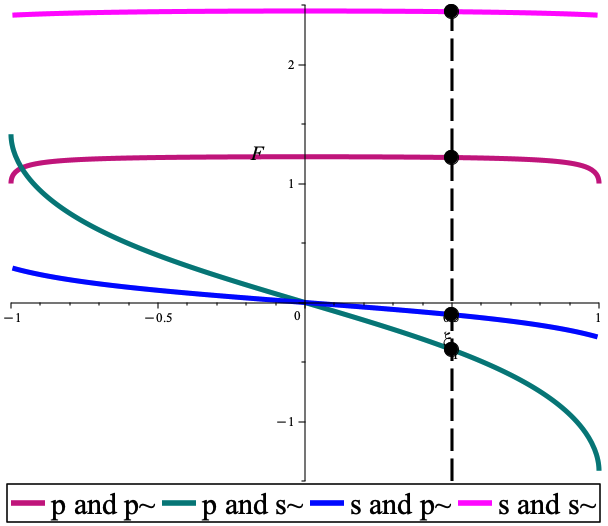}
        \subcaption{$F(\xi_1)$}
    \end{minipage}
    \caption{Example~\ref{ex:genpair2}(2) for $n=2$ and $\xi_0=1$. Left: the slowness surfaces of $(\rho,c)$, with pressure branch $p$ ($\abs{\xi}=1$) and shear branch $s$ ($\abs{\xi}=\sqrt2$), and of $(\tilde\rho,\tilde c)$, with branches $\tilde p$ ($\abs{\xi}=\sqrt{3/2}$) and $\tilde s$ ($\abs{\xi}=\sqrt3$). Choosing $\xi$ on $p$ or $s$ and $\tilde\xi$ on $\tilde p$ or $\tilde s$ with the same $\xi_1$, together with their polarizations ($\omega=\xi$ on $p$ and $\omega\perp\xi$ on $s$, likewise with tildes), gives a point of $W_m$, that is, $f=\tilde f=0$. The line $\xi_1=\frac12$ shows these choices for one value of $\xi_1$. As $\xi_1$ runs through $[-1,1]$, each of the four choices of branches traces an arc of one of the four irreducible components of $W_m$. Right: $F$ along these arcs as a function of $\xi_1$. Since $U_m=\{F=0\}$, a component lies in $U_m$ if and only if its curve vanishes identically. Isolated zeros, as at $\xi_1=0$, are allowed, but no curve vanishes identically, so no irreducible component of $W_m$ is contained in $U_m$.}
\label{fig:genpair2}
\end{figure}

The proof of Theorem~\ref{thm:main} relies on explicit examples of pairs $((\rho,c),(\tilde \rho, \tilde c))$ for each of the generic conditions in Lemma~\ref{lem:genPair} in general dimensions $n$ with $n\geq 2$. We illustrate concrete examples for both conditions in the following, the details can be found in the proofs of Lemma~\ref{lem:squarefreeness} and Lemma~\ref{lem:generic-noncontainment}:

\begin{Example}\label{ex:genpair2}
We give the following examples of generic pairs in Lemma~\ref{lem:genPair}:
\begin{enumerate}
\item[(1)] An example for Lemma~\ref{lem:genPair}(1): For positive $c$ and $\rho>0$, $(\rho,c)$ is squarefree if and only if $\Gamma_c(\xi)$ has $n$ distinct eigenvalues for some $\xi$. Consider
\begin{equation}
c^{(s)}_{ijk\ell}=\big(1+s(\min\{i,j\}-1)\big)\big(\delta_{ik}\delta_{j\ell}+\delta_{i\ell}\delta_{jk}\big),\qquad s\geq0,
\end{equation}
which interpolates between the isotropic tensor $c^{(0)}$ with $\lambda=0$, $\mu=1$ and $c^{(1)}=c^*$ from Lemma~\ref{lem:squarefreeness}; for $n=3$ its Voigt matrix is $\operatorname{diag}(2,\,2+2s,\,2+4s,\,1+s,\,1,\,1)$. As in the proof of Lemma~\ref{lem:squarefreeness}, $c^{(s)}$ satisfies~\eqref{eq:pos}, and
\begin{equation}
\Gamma_{c^{(s)}}(e_n)=\operatorname{diag}\big(1,\,1+s,\,\dots,\,1+(n-2)s,\,2+2(n-1)s\big)
\end{equation}
has $n$ distinct eigenvalues for $s>0$, so $(\rho,c^{(s)})$ is squarefree. For $s=0$ and $n\geq3$, the eigenvalue $\abs{\xi}^2$ of $\Gamma_{c^{(0)}}(\xi)=\xi\xi^T+\abs{\xi}^2I$ has multiplicity $n-1$ for every $\xi$, and $(\rho,c^{(0)})$ is not squarefree. Figure~\ref{fig:genpair1} illustrates this for $n=3$ and $\rho=1$: the $\xi_3$-axis meets the slowness surface in the six distinct points $\pm1$, $\pm(1+s)^{-1/2}$, $\pm(2+4s)^{-1/2}$ for $s>0$, whereas for $s=0$ the points $\pm1$ are double.
\item[(2)] An example for Lemma~\ref{lem:genPair}(2): Consider $\rho=2$, $\tilde\rho=3$ and the isotropic stiffness tensors $c_{ijk\ell}=\tilde c_{ijk\ell}=\delta_{ik}\delta_{j\ell}+\delta_{i\ell}\delta_{jk}$, i.e.\ $\lambda=\tilde\lambda=0$ and $\mu=\tilde\mu=1$. Then $\Gamma_c(\xi)=\xi\xi^T+(\xi\cdot\xi)I$ with $\xi\cdot\xi=\sum_{i=1}^n\xi_i^2$, so the slowness surface of $(\rho,c)$ consists of the pressure branch $p=\{\xi\cdot\xi=\xi_0^2\}$ with polarization $\omega=\xi$ and the shear branch $s=\{\xi\cdot\xi=2\xi_0^2\}$ with polarizations $\xi\cdot\omega=0$, and that of $(\tilde\rho,\tilde c)$ consists of $\tilde p=\{\xi\cdot\xi=\frac32\xi_0^2\}$ and $\tilde s=\{\xi\cdot\xi=3\xi_0^2\}$. Since $c=\tilde c$, we have $F=\xi_0^2\,\omega\cdot\tilde\omega$. As no branch of $(\rho,c)$ coincides with a branch of $(\tilde\rho,\tilde c)$, the irreducible components of $W_m$ are the four sets given by the choice of branches for $\xi$ and $\tilde\xi$. Let $\xi_0=\xi_1=1$, $\xi_2=\dots=\xi_{n-1}=0$, $\xi_n\in\{0,1\}$ and $\tilde\xi_n\in\{2^{-1/2},\sqrt2\}$, with $\omega=\xi$ on $p$, $\omega=\xi_ne_1-e_n$ on $S$, and likewise with tildes. This gives a point on each of the components $(p,\tilde p)$, $(p,\tilde s)$, $(s,\tilde p)$, $(s,\tilde s)$, with $F=1,\sqrt2,1-2^{-1/2},1+\sqrt2$, respectively. Hence no irreducible component of $W_m$ is contained in $U_m$. For $n=2$ the pair also satisfies condition~(1), since the two branches of each slowness surface are distinct; for $n\geq3$ the shear branch occurs $n-1$ times and condition~(1) fails. Figure~\ref{fig:genpair2} shows the case $n=2$, where each $\xi_1\in[-1,1]$ gives a point with $\xi_2,\tilde\xi_2>0$ on each of the four branches.
\end{enumerate}
   
\end{Example}

We use a construction of special solutions as presented in Proposition~\ref{prop_concentrating_sol} to obtain the following Lemma by using the Alessandrini identity for this construction. This Lemma is a reformulated version of Proposition~\ref{prop:Alessandrini}:
\begin{lemma}\label{lem:Alessandrini}
    Let $\hat{x}=(\hat{x}',0)\in (0,T)\times \partial \Omega$, where $(x_0, x_1, \ldots, x_n)$ are coordinates in $(0,T) \times \overline{\Omega}$ with $x_n = 0$ corresponding to $(0,T) \times \partial \Omega$ near $\hat{x}$, and denote by $((\rho,c),(\tilde \rho, \tilde c))$ a pair of stiffness tensors and densities. Suppose that there is $\hat{\xi}=(\hat{\xi}',0)$ such that $\spoly{\rho,c}(\hat{x}',0,\hat{\xi}',z)$ and $\spoly{\tilde \rho, \tilde c}(\hat{x}',0,\hat{\xi}',\tilde z)$ only have simple non-real roots in $z$ and $\tilde z$ respectively (coming in conjugate pairs)
. Let $\omega$ and $\tilde \omega$ be the unique eigenvectors such that $f_i(\hat{\xi}',z,\hat{\xi}',\tilde z, \omega, \tilde \omega)=\tilde f_i(\hat{\xi}',z,\hat{\xi}',\tilde z, \omega, \tilde \omega)=0$ for all $i=1,\dots, n$. If $\Lambda_{\rho,c}=\Lambda_{\tilde \rho, \tilde c}$ then
    \begin{equation}
((c(x',0)-\tilde c(x',0))_{ijkl} \zeta_k \overline{\tilde{\zeta}_j} - (\rho(x',0)-\tilde \rho(x',0))\delta_{il} \xi_0^2) \omega_l \overline{\tilde{\omega}_i} = 0
\end{equation}
for $(x',\xi')$ near $(\hat{x}', \hat{\xi}')$. Here 
\begin{equation}
\zeta(x',\xi') = (\xi', z(x',\xi')), \qquad \tilde{\zeta}(x',\xi') = (\xi', \tilde{z}(x',\xi')).
\end{equation}
\end{lemma}

Since $\hat\xi'$ is real and $\tilde c$, $\tilde\rho$ are real, $\overline{\tilde z}$ is again a root of $\spoly{\tilde\rho,\tilde c}(\hat x',0,\hat\xi',\,\cdot\,)$, with kernel vector $\overline{\tilde\omega}$, and $\overline{\tilde\zeta}=(\xi',\overline{\tilde z})$. Hence, at $x'=\hat x'$, the identity in Lemma~\ref{lem:Alessandrini} states that $F$ from~\eqref{eq:Fintro} vanishes at the point $(\xi,\tilde\xi,\omega,\tilde\omega)=(\zeta,\overline{\tilde\zeta},\omega,\overline{\tilde\omega})$, at which also $f_i=\tilde f_i=0$. The only information passed from the analysis to the algebraic part is that $F$ vanishes at such points. Therefore $F$, $f_i$ and $\tilde f_i$ can be defined without complex conjugates; note that $\tilde\xi_n=\overline{\tilde z}$ has negative imaginary part. 

The setting in this Lemma motivates the set $\Xi$ defined in~\eqref{eq:Xidef}, as $\zeta$ and $\tilde \zeta$ have the same entries for $i=1,\dots,n-1$. Additionally this Lemma motivates to define the variety $W_m$ in~\eqref{eq:Wmeq} and the variety $U_m$ in~\eqref{eq:Umeq} to capture those $\xi,\tilde \xi, \omega, \tilde \omega\in \Xi \times (\C\P)^2$ for which $f_i$ and $\tilde f_i$ vanish for all $i=1,\dots, n$ and $F$ vanishes respectively.\par
With these Lemmas we are ready to give the following sketch of a proof for Theorem~\ref{thm:main} and we refer the reader to Theorem~\ref{thm:generic-uniqueness} for the full proof:

\begin{proof}[Sketch of the proof of Theorem~\ref{thm:main}]
Consider a generic pair $((\rho,c),(\tilde \rho, \tilde c))$ as in Lemma~\ref{lem:genPair}.
When $(\rho,c)$ and $(\tilde \rho, \tilde c)$ individually give rise to a square-free slowness polynomial this is equivalent for them to satisfy an explicit generic condition on the pair (see Lemma~\ref{lem:generic-simplicity}(2)). The slowness polynomial $\spoly{\rho,c}(\xi',z,\xi_0)$ of a generic stiffness tensor $c$ and density $\rho$ is a square-free polynomial in $\C[\xi_0,\xi',z]$ if and only if $\spoly{\rho,c}$ is square-free in $\C[z]$ (see Lemma~\ref{lem:squareFree}). The latter condition is equivalent to $\spoly{\rho,c}(\xi',z,\xi_0)$ having simple roots in $z$. This implies for each entry of the pair that $\spoly{\rho, c}(x,\xi',z)$ and $\spoly{\tilde \rho, \tilde c}(x,\xi',\tilde z)$ have simple roots in $z$ and $\tilde z$ respectively. Because the stiffness tensors $c$ and $\tilde c$ are positive definite, the real slowness surface is compact such that if one chooses $(\hat\xi_1,\dots,\hat\xi_{n-1})$ large enough compared to~$\hat\xi_0$ the roots $z$ and $\tilde z$ of $\spoly{\rho, c}(x,\hat{\xi}',z)$ and $\spoly{\tilde \rho, \tilde c}(x,\hat{\xi}',\tilde z)$ will be non-real. As shown in the proof of Lemma~\ref{lem:alg-alessandrini}, $\hat\xi'$ can be chosen such that in addition all roots are simple. Hence, by the generic condition on square-freeness in Lemma~\ref{lem:genPair} we can apply Lemma~\ref{lem:Alessandrini} that gives the following identity for each boundary point $\hat{x}=(\hat{x}',0)\in (0,T)\times \partial \Omega$ 
\begin{equation}
F(\zeta,\overline{\tilde\zeta},\omega,\overline{\tilde\omega}) = 0
\end{equation}
for $\xi'$ near $\hat{\xi}'$, with $F$ as in~\eqref{eq:Fintro}, $\zeta$, $\tilde\zeta$ as in Lemma~\ref{lem:Alessandrini}, and for $\omega$ and $\tilde \omega$ being the unique eigenvectors satisfying
\begin{equation}
f_i(\hat{\xi}',z,\hat{\xi}',\tilde z, \omega, \tilde \omega)=\tilde f_i(\hat{\xi}',z,\hat{\xi}',\tilde z, \omega, \tilde \omega)=0
\end{equation}
for all $i=1,\dots, n$.
It does not follow for general anisotropic stiffness tensors $c$ and $\tilde c$ and densities $\rho$ and $\tilde \rho$ that $F=0$ implies $c=\tilde c$ and $\rho =\tilde \rho$, therefore we use additional tools from algebraic geometry to obtain the desired result. In the setting of Lemma~\ref{lem:Alessandrini} all entries of $\xi$ and $\tilde \xi$ are real apart from the $n$th entry; for that purpose we introduce the following subset:
\begin{equation}
A_\R
=
\{
(\xi,\tilde\xi,\omega,\tilde\omega)
\in A
;
\xi_i\in\R
\text{ for all }
i=0,\dots,n-1
\}
.
\end{equation}
We now consider a point $p=(\xi, \tilde \xi, \omega, \tilde \omega)$ in a non-empty Euclidean open subset $O \subset A_{\R}\cap W_m$ in the Euclidean topology so that also $O\subset U_m$. Since the point $p=(\xi,\tilde \xi, \omega, \tilde \omega)\in W_m$ that corresponds to a squarefree pair $((\rho,c),(\tilde \rho, \tilde c))$ is a smooth point of $W_m$ (see Lemma~\ref{lem:smooth}), the irreducible components of $W_m$ are disjoint. Let $C$ be the irreducible component containing~$p$. $O\cap C$ contains a non-empty Euclidean open subset and this subset is Zariski dense in $C$. Since $O\subset U_m$ and $U_m$ is Zariski closed it follows that $C\subset U_m$, contradicting Lemma~\ref{lem:genPair}(2) that no irreducible component of~$W_m$ is contained in~$U_m$.
\end{proof}




The corollary for fields instead of point values is proven in section~\ref{sec:fields}, and it is a standard topological argument built on the pointwise result.

\section{Construction of special solutions}
\label{sec:specialSols}

In this section, we will construct concentrating solutions to extract boundary information on coefficients of a system from boundary measurements of solutions. The main boundary determination result from the Cauchy data set is Theorem \ref{thm_bd}, which extends a result in \cite{Oksanen2024} from scalar equations to systems. As this result holds for more general operators than the elastic wave operator, we phrase everything in a more general setting in this section. Later in Proposition \ref{prop:Alessandrini} we will give a boundary determination result from the DN map of the elastic wave equation.

Let $M$ be a compact $n$-manifold with smooth boundary, and let $P$ be an $N \times N$ matrix whose entries are differential operators (with smooth coefficients) of order $m$ on $M$. We wish to study boundary determination of coefficients of $P$ from boundary measurements. Let $\hat{x} \in \p M$ be a boundary point. Boundary determination in the elliptic region as in \cite[Section 5]{Oksanen2024} involves the construction of solutions of $Pu = 0$ in $M$ that concentrate near $\hat{x}$ and oscillate rapidly on the boundary. Let us describe this construction.

We use a system of coordinates $(x',x_n)$ near $\hat{x}$ so that $\hat{x} = 0$ and near $\hat{x}$ the set $M$ is given by $\{ x_n \geq 0 \}$. In these coordinates the operator $P$ may be written as 
\begin{equation}
(P u)_r = P_{rs}(x,D) u_s, \qquad P_{rs}(x,D) = \sum_{|\alpha| \leq m} p_{rs,\alpha}(x) D^{\alpha}.
\end{equation}
For $h > 0$ small, we look for an approximate solution (satisfying $Pu = O(h^{\infty})$) in the form 
\begin{equation}
u = e^{i\Phi/h} a,
\end{equation}
where $\Phi$ is a smooth complex scalar function and $a$ is a smooth complex $N$-vector.

We have the conjugation formula \cite[Lemma 3.3]{Oksanen2024} 
\begin{equation} \label{conjugation_scalar}
e^{-i\Phi/h} P_{rs} (e^{i\Phi/h} a_s) = h^{-m} p_{rs,m}(x,\nabla \Phi) a_s + h^{1-m}(\frac{1}{i} L_{rs} a_s + b_{rs} a_s) + \sum_{k=2}^m h^{k-m} R_k a
\end{equation}
where $p_{rs,m}$ is the principal symbol of $P_{rs}$, $R_k$ is a matrix differential operator of order $k$, and 
\begin{align}
L_{rs} &= \p_{\xi_j} p_{rs,m}(x,\nabla \Phi) \p_j, \\
b_{rs} &= \frac{1}{2i} \p_{\xi_j \xi_k} p_{rs,m}(x,\nabla \Phi) \p_{jk} \Phi + p_{rs,m-1}(x,\nabla \Phi)
\end{align}
where the full symbol of $P_{rs}$ in these coordinates has the expansion $p_{rs} = p_{rs,m} + p_{rs,m-1} + \ldots$ with $p_{rs,k}(x,\xi)$ being homogeneous of degree~$k$ in~$\xi$.

We decompose~$a$ as $a = a_0 + h a_1 + h^2 a_2 + \ldots$ and write $a_j = (a_{j,1}, \ldots, a_{j,N})$. Grouping like powers of $h$ in~\eqref{conjugation_scalar}, we would like to solve the following equations to infinite order at $\{x_n=0\}$:
\begin{align} \label{eikonal_transport_eq}
\left\{ \begin{array}{rl}
p_m(x,\nabla \Phi) a_0 &\!\!\!= 0, \\[3pt]
\frac{1}{i} L_{rs} a_{0,s} + b_{rs} a_{0,s} + (p_m(x,\nabla \Phi) a_1)_r &\!\!\!= 0, \\[3pt]
\frac{1}{i} L_{rs} a_{1,s} + b_{rs} a_{1,s} + (p_m(x,\nabla \Phi) a_2)_r + R_2 a_0 &\!\!\!= 0,\\[3pt]
 &\hspace{-3pt}\vdots
 \end{array} \right.
\end{align}

Let us give some conditions for $P$ so that we can solve these equations. Recall that our boundary point of interest is $\hat{x} = (\hat{x}',0) = 0$. Covectors in $T_{\hat{x}}(\p M)$ can be written as $(\xi',0)$ or just as $\xi'$. We consider the set 
\begin{equation}
\Xi(\hat{x}) = \Xi(P, \hat{x})
= \{ (\xi',z) \in (T_{\hat{x}}(\p M) \setminus 0) \times \mC \,:\,  \det(p_m(\hat{x}',0,\xi',z)) = 0 \text{ and } \im(z) > 0 \}.
\end{equation}
That is, we wish to consider $\xi'$ for which the map 
\begin{equation}
z \mapsto \det(p_m(\hat{x}',0,\xi',z))
\end{equation}
has a complex root $z$ with $\im(z) > 0$. Since the solutions constructed in \cite{Oksanen2024} look roughly like $e^{i(x' \cdot \xi' + z x_n + O(x_n^2))/h} a$, the pairs $(\xi',z) \in \Xi(\hat{x})$ are the ones for which one might be able to construct solutions concentrating at $\hat{x}$, oscillating in direction $\xi'$ on $\p M$, and decaying exponentially in the interior.

To perform the asymptotic construction successfully, we further restrict to pairs $(\xi',z)$ such that the root $z$ of $\det(p_m(\hat{x}',0,\xi',\,\cdot\,))$ is simple. Consider
\begin{equation}
\Xi_{1}(\hat{x}) = \Xi_{1}(P,\hat{x})
 = \{ (\xi', z) \in \Xi(\hat{x}) \,:\,  \text{$z$ is a simple root of $\det(p_m(\hat{x}',0,\xi',\,\cdot\,))$} \}.
\end{equation}
Note that the simple root condition is equivalent with
\begin{equation}
\p_w (\det(p_m(\hat{x}',0,\xi',\,\cdot\,)))|_{w=z} \neq 0,
\end{equation}
and that, by Jacobi's formula, it implies $\dim(\ker(p_m(\hat{x}',0,\xi',z))) = 1$, see the proof of Lemma~\ref{lemma_z_omega}.

Next we show some properties related to pairs $(\hat{\xi}',\hat{z}) \in \Xi_{1}(\hat{x})$ that are valid when varying $(\hat{x}', \hat{\xi}')$ slightly.

\begin{lemma} \label{lemma_z_omega}
Suppose that $(\hat{\xi}',\hat{z}) \in \Xi_{1}(\hat{x})$. For any $x'$ near $\hat{x}'$ and any $\xi'$ near $\hat{\xi}'$ there is a pair $(\xi',z(x',\xi'))$ in $\Xi_{1}(x',0)$ and a unit vector $\omega(x',\xi') \in \mC^N$ such that 
\begin{equation}
\ker(p_m({x}',0,\xi',z(x',\xi'))) = \mC \omega(x',\xi').
\end{equation}
One can choose $z$ and $\omega$ depending smoothly on $(x',\xi')$. Moreover, $\operatorname{Ran}(p_m(\,\cdot\,)) = \gamma^{\perp}$ for some unit vector $\gamma = \gamma(x',\xi') \in \mC^N$ depending smoothly on $(x',\xi')$, and  one has 
\begin{equation}
(\p_{\xi_n} p_m(\,\cdot\,) \omega, \gamma) \neq 0.
\end{equation}
Given a vector $b \in \gamma^{\perp}$ depending smoothly on $(x',\xi')$, there is a unique vector $y$ solving 
\begin{equation}
p_m(\,\cdot\,) y = b, \qquad y \perp \omega,
\end{equation}
and $y$ depends smoothly on $(x',\xi')$.
\end{lemma}
\begin{proof}
The fact that $z$ and $\omega$ exist and depend smoothly on $(x',\xi')$ is a standard consequence of the implicit function theorem and the fact that the roots are assumed to be simple. We give the details for completeness. First consider the function 
\begin{equation}
F(x',\xi',z) = \det(p_m(x',0,\xi',z)).
\end{equation}
This function is smooth near $(\hat{x}',\hat{\xi}',\hat{z})$, satisfies $F(\hat{x}',\hat{\xi}',\hat{z}) = 0$ and $\p_z F(\hat{x}',\hat{\xi}',\hat{z}) \neq 0$. By the implicit function theorem (with real parameters $x',\xi'$ and holomorphic parameter $z$), there is a smooth function $z = z(x',\xi')$ such that $F(x',\xi',z(x',\xi'))=0$ near $(\hat{x}',\hat{\xi}')$ and $z(\hat{x}',\hat{\xi}') = \hat{z}$. For $(x',\xi')$ near $(\hat{x}',\hat{\xi}')$, by continuity $\p_w F({x}',{\xi}',\,\cdot\,)|_{w=z(x',\xi')} \neq 0$ so $z(x',\xi')$ is a simple root, and $\im(z(x',\xi')) > 0$, so that $(\xi',z(x',\xi')) \in \Xi_1(x',0)$.

Write $A = A(x',\xi') = p_m({x}',0,\xi',z(x',\xi'))$. The simple root condition and Jacobi's formula imply that
\begin{align}
0 &\neq \p_w (\det(p_m({x}',0,\xi',\,\cdot\,)))|_{w=z} = \operatorname{tr}(\p_w(p_m({x}',0,\xi',\,\cdot\,))|_{w=z} \operatorname{adj}(A)) \label{simple_root_jacobi} \\
 &=  \operatorname{tr}(\p_{\xi_n} p_m({x}',0,\xi',z) \operatorname{adj}(A)). \notag
\end{align}
In the last equality we used analyticity in $w$. In particular $\operatorname{adj}(A) \neq 0$, and since $\det(A) = 0$ this gives $\operatorname{rank}(A) = N-1$. Thus $\ker(A)$ is one-dimensional, $\operatorname{Ran}(A)$ is $(N-1)$-dimensional, and $\operatorname{adj}(A)$ has rank one. From $A \operatorname{adj}(A) = \operatorname{adj}(A) A = \det(A) \operatorname{Id} = 0$ we see that every non-zero column of $\operatorname{adj}(A)$ spans $\ker(A)$ and every non-zero row of $\operatorname{adj}(A)$ spans $\ker(A^T)$. Choose a column of $\operatorname{adj}(\hat{A})$, $\hat{A} = p_m(\hat{x}',0,\hat{\xi}',\hat{z})$, that is non-zero. By continuity the corresponding column of $\operatorname{adj}(A(x',\xi'))$ is non-zero for $(x',\xi')$ near $(\hat{x}',\hat{\xi}')$, and its entries are polynomials in the entries of $A(x',\xi')$ and hence smooth. Normalizing it gives a smooth unit vector $\omega(x',\xi')$ with $\ker(A) = \mC \omega$. In the same way, normalizing the complex conjugate of a non-zero row of $\operatorname{adj}(A)$ gives a smooth unit vector $\gamma(x',\xi')$ with $A^* \gamma = 0$, and since $\operatorname{Ran}(A) = \ker(A^*)^{\perp}$ is $(N-1)$-dimensional, we have $\operatorname{Ran}(A) = \gamma^{\perp}$. As the columns of the rank one matrix $\operatorname{adj}(A)$ are multiples of $\omega$ and its rows are multiples of $\bar{\gamma}^T$, we may write
\begin{equation}
\operatorname{adj}(A) = c \omega \bar{\gamma}^T
\end{equation}
for some $c = c(x',\xi') \in \mC$. By~\eqref{simple_root_jacobi} we have $c \neq 0$. Inserting this in~\eqref{simple_root_jacobi} gives
\begin{equation}
0 \neq c \,\operatorname{tr}(\p_{\xi_n} p_m(\,\cdot\,) \omega \bar{\gamma}^T) = c (\p_{\xi_n} p_m(\,\cdot\,) \omega, \gamma).
\end{equation}
This proves that $(\p_{\xi_n} p_m(\,\cdot\,) \omega, \gamma) \neq 0$.

Finally, let $b \in \gamma^{\perp}$ depend smoothly on $(x',\xi')$ and consider the map 
\begin{equation}
H(x',\xi',y) = Ay - b + (y \cdot \omega) \gamma.
\end{equation}
This map is smooth. Since $\operatorname{Ran}(A) = \gamma^{\perp}$, for any $(x',\xi')$ there is a unique vector $y(x',\xi')$ satisfying $H(x',\xi',y)=0$, and $y \in \omega^{\perp}$. Writing $\hat{y} = y(\hat{x}',\hat{\xi}')$, we have 
\begin{equation}
D_y H(\hat{x}',\hat{\xi}', \hat{y}) v = Av + (v \cdot \omega)\gamma|_{(\hat{x}', \hat{\xi}')}.
\end{equation}
This vanishes if and only if $v=0$. Thus $D_y H$ is invertible at this point, and the implicit function theorem ensures that $y(x',\xi')$ is smooth near $(\hat{x}',\hat{\xi}')$.
\end{proof}
 
Next we show that for any $(\hat{\xi}',\hat{z}) \in \Xi_1$ (if such a pair exists), one can construct concentrating solutions as in \cite[Theorem 5.1]{Oksanen2024}.

\begin{proposition} \label{prop_concentrating_sol}
Let $P$ be a $N \times N$ matrix of differential operators of order $m \geq 1$. Let $\hat{x} = 0 \in \p M$ in coordinates as above, and let $(\hat{\xi}',\hat{z}) \in \Xi_{1}(\hat{x})$. There exists $u = e^{i\Phi/h} a$ where $\Phi \in C^{\infty}(M,\mC)$, $a \sim \sum h^k a_k$ with $a, a_k \in C^{\infty}(M,\mC^N)$ and $\supp(a)$ near $\hat{x}$, such that the following conditions hold for $(x',\xi')$ near $(\hat{x}',\hat{\xi}')$ with $\Phi, a_j$ depending also on $\xi'$:
\begin{align}
\Phi(x',0) &= x' \cdot \xi', \\
\p_n \Phi(x',0) &= z(x',\xi') , \\
a_0(x',0) &= \psi(x') \omega(x',\xi'), \\
\im(\Phi) &\geq 0 \text{ on $\supp(a)$},
\end{align}
where $\psi(x')$ is a $C^{\infty}$ function with small support such that $\psi = 1$ near $\hat{x}'$. Moreover, if $s \geq 0$ is fixed, one has 
\begin{equation}
\norm{u}_{L^2(M)} \sim h^{1/2}, \qquad \norm{u}_{H^k(M)} \leq C_k h^{-k+1/2}, \qquad \norm{Pu}_{H^s(M)} = O(h^{\infty}).
\end{equation}
\end{proposition}
\begin{proof}
Let $(\xi', z(x',\xi')) \in \Xi_{1}(x',0)$ as in Lemma \ref{lemma_z_omega}. From~\eqref{eikonal_transport_eq}, we first need to find $\Phi$ and $a_0$ such that  
\begin{equation} \label{eikonal_inf}
p_m(x,\nabla \Phi) a_0 = 0 \text{ to infinite order on $\{x_n=0\}$}.
\end{equation}
We first require that 
\begin{align}
\Phi(x',0) &= x' \cdot \xi', \\
\p_n \Phi(x',0) &= z(x',\xi').
\end{align}
Then $\ker(p_m(x',0,\nabla \Phi(x',0))) = \mC \omega(x',\xi')$, and we initially choose 
\begin{equation}
a_0(x',0) = \omega(x',\xi').
\end{equation}
It follows that~\eqref{eikonal_inf} is satisfied to order $0$.

We proceed by applying $\p_{x_n}$ to the expression $p_m(x,\nabla \Phi) a_0$  and evaluating at $x_n=0$. This gives the equation  
\begin{equation} \label{eikonal_first_order}
p_m(\,\cdot\,) (\p_n a_0) + \p_{x_n} p_m(\,\cdot\,) a_0 + \p_{\xi_{\alpha}} p_m(\,\cdot\,) (\p_{\alpha n} \Phi) a_0 +  \p_{\xi_{n}} p_m(\,\cdot\,) (\p_{n n} \Phi) a_0|_{x_n=0} = 0.
\end{equation}
Write $A = A(x',\xi') = p_m(\,\cdot\,)|_{x_n=0} = p_m(x',0,\xi',z(x',\xi'))$. As in Lemma \ref{lemma_z_omega} we have $\operatorname{Ran}(A) = \gamma^{\perp}$ for some unit vector $\gamma$ with 
\begin{equation}
(\p_{\xi_n} p_m(\,\cdot\,) \omega, \gamma) \neq 0.
\end{equation}
Taking the inner product of~\eqref{eikonal_first_order} with $\gamma$, we will first choose $\p_{nn} \Phi(x',0)$ such that 
\begin{equation}
(\p_{\xi_{n}} p_m(\,\cdot\,) \omega, \gamma) \p_{n n} \Phi(x',0) = -(\p_{x_n} p_m(\,\cdot\,) \omega + \p_{\xi_{\alpha}} p_m(\,\cdot\,) (\p_{\alpha n} \Phi) \omega, \gamma).
\end{equation}
Since we may divide by $(\p_{\xi_{n}} p_m(\,\cdot\,) \omega, \gamma) \neq 0$, it is possible to find $\p_{n n} \Phi(x',0)$ with this property. On the other hand, since $\operatorname{Ran}(A) = \gamma^{\perp}$, we can arrange $p_m(\,\cdot\,) (\p_n a_0)|_{x_n=0}$ to be any vector in $\gamma^{\perp}$ by a suitable choice of $\p_n a_0|_{x_n=0}$. Furthermore, the last part of Lemma \ref{lemma_z_omega} guarantees that $\p_n a_0|_{x_n=0}$ can be chosen to depend smoothly on $(x',\xi')$. With these choices~\eqref{eikonal_inf} holds to order $1$. Iterating this construction allows us to find the Taylor series of $\Phi$ and $a_0$ on $\{x_n=0\}$ so that~\eqref{eikonal_inf} holds to infinite order. Then Borel summation gives the required functions $\Phi$ and $a_0$.

We now move to the construction of $a_1$. Note that if~\eqref{eikonal_inf} holds, then it also holds with $a_0$ replaced by $c_0 a_0$ where $c_0$ is a smooth scalar function. This gives us more flexibility in the construction. The first transport equation in~\eqref{eikonal_transport_eq} becomes 
\begin{equation} \label{transport_first}
\frac{1}{i} L_{rs} (c_0 a_{0,s}) + b_{rs} c_0 a_{0,s} + (p_m(x,\nabla \Phi) a_1)_r = 0.
\end{equation}
We wish to find the Taylor series of $c_0$ and $a_1$ at $x_n=0$ so that~\eqref{transport_first} holds  to infinite order at $x_n=0$. We first solve~\eqref{transport_first} to order $0$ at $x_n=0$. As above, we have $\operatorname{Ran}(p_m(\,\cdot\,)|_{x_n=0}) = \gamma^{\perp}$ and $a_0(x',0) = \omega$. Thus by making a suitable smooth choice of $a_1(x',0)$ using the last part of Lemma \ref{lemma_z_omega}, it is enough to find $c_0$ satisfying 
\begin{equation}
(\frac{1}{i} L_{rs} c_0) \omega_s + (b_{rs} \omega_s + \frac{1}{i} L_{rs} \omega_s)c_0) \bar{\gamma}_r \big|_{x_n=0}= 0.
\end{equation}
This can be written as the scalar equation $\frac{1}{i} Lc_0 + q c_0|_{x_n=0} = 0$, where 
\begin{equation}
L = \omega_s \bar{\gamma}_r L_{rs} = (\p_{\xi_j} p_{rs,m}\omega_s \bar{\gamma}_r ) \p_j.
\end{equation}
Since $(\p_{\xi_n} p_m(\,\cdot\,) \omega, \gamma)|_{x_n=0} \neq 0$ by Lemma \ref{lemma_z_omega}, we can prescribe $c_0(x',0) = \psi(x')$ and solve for $\p_n c_0(x',0)$. Note that $c_0 a_0|_{x_n=0} = \psi(x') \omega(x',\xi')$ as required. With these choices, the equation~\eqref{transport_first} holds to order $0$. Applying $\p_{x_n}$ to~\eqref{transport_first} and evaluating at $x_n=0$ leads to an equation of the same form involving $\p_n^2 c_0(x',0)$ and $\p_n a_1(x',0)$, and this equation can be solved by the same procedure as above. Iterating this construction and using Borel summation gives us smooth functions $c_0$ and $a_1$ such that~\eqref{transport_first} holds to infinite order at $x_n=0$. The amplitudes $a_2, a_3, \ldots$ can be constructed in the same way.

We have now constructed $u = e^{i\Phi/h} a$ such that $Pu = O(x_n^{\infty})$ near~$\hat{x}$. If we replace~$a$ by~$\chi a$ where $\chi(x',x_n) = \chi_1(x') \chi_2(x_n)$ is a smooth cutoff function with $\chi_1 = 1$ near $\supp(\psi)$, $\chi_2 = 1$ near~$0$, and $\operatorname{Im}(\Phi) \geq 0$ in~$\supp(\chi)$, then~$u$ will satisfy the required properties.
\end{proof}

\begin{remark}
It is possible that the construction would work also for non-square systems~$P$ and for pairs~$\xi'$ and~$z$ with $\im(z) > 0$ under a suitable simplicity condition and assuming that 
\begin{equation}
(\p_{\xi_n} p_m(\,\cdot\,) \omega, \gamma) \neq 0
\end{equation}
where $\ker(p_m(\,\cdot\,)) = \mC \omega$ and $\operatorname{Ran}(p_m(\,\cdot\,)) = \gamma^{\perp}$.
\end{remark}

\section{Boundary determination from Cauchy data set}\label{sec:bdCauchy}

We now show how one can use the concentrating solutions given by Proposition \ref{prop_concentrating_sol} to prove boundary determination from Cauchy data for certain matrix elements of the principal symbol. If $P$ is an $N \times N$ matrix of differential operators in $M$ of order $m$, we define the Cauchy data set 
\begin{equation}
C_P = \{ (u|_{\p M}, \nabla u|_{\p M}, \ldots, \nabla^{m-1} u|_{\p M}) \,:\, u \in H^m(M, \C^N) \text{ satisfies } Pu = 0 \text{ in $M$} \}.
\end{equation}
Here~$\nabla$ is the total covariant derivative corresponding to some fixed Riemannian metric on~$M$. We assume that one can correct approximate solutions to exact ones, in the sense that there exist $s \in \R$ and $C > 0$ such that for any $u_0 \in H^{m+s}(M, \C^N)$, the equation 
\begin{equation}
Pv = -Pu_0 \quad \text{in $M$}
\end{equation}
has a solution $v \in H^m(M)$ with 
\begin{equation}
\norm{v}_{H^m(M)} \leq C \norm{Pu_0}_{H^s(M)}.
\end{equation}
Such an estimate holds e.g.\ for the fully anisotropic elastic wave operator in a space-time cylinder~\cite{stolk-thesis}.

Our first theorem considers boundary determination for certain matrix elements of the principal symbol of a PDE system.


\begin{theorem} \label{thm_bd}
Let $P$ and $Q$ be two $N \times N$ matrices of differential operators of order $m \geq 1$. Let $\hat{x} = (\hat{x}',0) \in \p M$, and suppose that there is $\hat{\xi}' = (\hat{\xi}',0) \in T_{\hat{x}}(\p M) \setminus 0$ such that 
\begin{equation}
(\xi', z) \in \Xi_1(P, (x',0)), \qquad (\xi',\rho) \in \Xi_{1}(Q^*, (x',0)),
\end{equation}
where $z$ and $\rho$ depend smoothly on $(x',\xi')$ near $(\hat{x}',\hat{\xi}')$. Let $\ker(p_m(x',0,\xi',z)) = \mC \omega$ and $\ker(q_m^*(x',0,\xi',\rho)) = \mC \beta$. If 
\begin{equation}
C_P = C_Q
\end{equation}
and if $P$ and $Q^*$ are such that one can correct approximate solutions to exact ones as described above, then 
\begin{gather*}
p_m(x',0,\xi',z) \omega = 0,  \qquad 
((p_m - q_m)(x',0,\xi',z) \omega, \beta) = 0,
\end{gather*}
for $(x',\xi')$ near $(\hat{x}', \hat{\xi}')$.
\end{theorem}
\begin{proof}
The proof is analogous to the corresponding result in \cite{Oksanen2024}. As in \cite[Lemma 2.8]{Oksanen2024}, the condition $C_P = C_Q$ implies the integral identity 
\begin{equation} \label{bd_integral_identity}
((P-Q) u_1, u_2)_{L^2(M)} = 0
\end{equation}
for any~$H^m$ solutions of $Pu_1 = 0$ and $Q^* u_2 = 0$ in~$M$.

By Proposition  \ref{prop_concentrating_sol} we construct approximate solutions 
\begin{equation}
v_1 = h^{-1/2} e^{i\Phi_1/h} a_1, \qquad v_2 = h^{-1/2} e^{i\Phi_2/h} a_2
\end{equation}
for the operators $P$ and $Q^*$ corresponding to $(\xi', z) \in \Xi_1(P, (x',0))$ and $(\xi',\rho) \in \Xi_{1}(Q^*, (x',0))$, respectively. The functions $v_j$ satisfy $\norm{v_j}_{L^2} \sim 1$, $\|Pv_1\|_{H^s} = O(h^{\infty})$, $\|Q^* v_2\|_{H^s} = O(h^{\infty})$ with $\Phi_1(x', 0) = x' \cdot \xi'$, $\p_n \Phi_1(x',0) = z(x',\xi')$, $a_{1,0} = \psi_1(x') \omega(x',\xi')$, and $\im(\Phi_1) \geq 0$ in $\supp(a_{1})$, where $\psi_1$ is a smooth real function near $\hat{x}'$. The function $v_2$ has similar properties but with $z$ replaced by $\rho$, $\omega$ replaced by $\beta$, and $\psi_1$ replaced by $\psi_2$. By the well-posedness assumption in the lemma, we obtain exact solutions 
\begin{equation}
u_j = v_j + r_j
\end{equation}
where $\norm{r_j}_{H^m(M)} = O(h^{\infty})$.

Inserting the solutions~$u_j$ in~\eqref{bd_integral_identity} and applying the conjugation formula~\eqref{conjugation_scalar} implies that 
\begin{align}
0 &= \lim_{h \to 0} h^m \int_M (P-Q)u_1 \cdot \bar{u}_2 \,dV \\
 &=\lim_{h \to 0} h^{-1} \int_M (p_{m}(x,d\Phi_1)-q_{m}(x,d\Phi_1))e^{i(\Phi_1-\bar{\Phi}_2)/h} a_1 \cdot \bar{a}_2 \,dV.
\end{align}
Write $r = p_{m}(x,d\Phi_1)-q_{m}(x,d\Phi_1)$ and $\Theta = \Phi_1-\bar{\Phi}_2$. In boundary normal coordinates near~$\hat{x}$ we have 
\begin{equation}
\Theta = (z-\bar{\rho}) x_n + O(x_n^2)
\end{equation}
where $\im(z-\bar{\rho}) > 0$. If the supports of $a_j$ are chosen small enough, it follows that 
\begin{align}
0 &= \lim_{h \to 0} \, h^{-1} \int_M r e^{i\Theta/h} a_1 \cdot \bar{a}_2 \,dV \\
 &= \lim_{h \to 0} \, h^{-1} \int_{\{ x_n > 0 \}} e^{i((z-\bar{\rho}) x_n + O(x_n^2))/h} r a_1 \cdot \bar{a}_2 |g|^{1/2} \,dx \\
 &= \lim_{h \to 0} \, \int_{\{ x_n > 0 \}} e^{i((z-\bar{\rho}) x_n + h O(x_n^2))} (r a_1 \cdot \bar{a}_2 |g|^{1/2})(x',hx_n) \,dx \\
 &= \int_{\{ x_n > 0 \}} e^{i((z-\bar{\rho}) x_n)} (r a_{1,0} \cdot \bar{a}_{2,0} |g|^{1/2})(x',0) \,dx \\
 &= \int_{\R^{n-1}} \frac{i}{z-\bar{\rho}} \psi_1 \psi_2 |g|^{1/2} (r\omega \cdot \bar{\beta})(x',0) \,dx'.
\end{align}
Since $z-\bar{\rho}$ and $|g|^{1/2}$ are nonvanishing near $0$, varying the test functions $\psi_1$ and $\psi_2$ implies that $(r \omega \cdot \bar{\beta})(x',0) = 0$. Using the properties of $\Phi_1$, this can be rewritten as 
\begin{equation}
(p_m-q_m)(x',0,\xi',z(x',\xi')) \omega(x',\xi') \cdot \bar{\beta}(x',\xi') = 0
\end{equation}
for $(x',\xi')$ near $(\hat{x}', \hat{\xi}')$. This proves the claim.
\end{proof}

Theorem \ref{thm_bd} states that if $P$ and $Q$ have the same boundary measurements, and if appropriate nondegeneracy conditions hold, then 
\begin{equation}
(p_m - q_m) \omega \cdot \beta = 0 \text{ on $\partial M$}
\end{equation}
for certain vectors $\omega$ and $\beta$ (and related roots $z$ and $\rho$). If we are in a setting where one has many linearly independent vectors and distinct roots, then it is possible to determine the characteristic polynomial of the principal symbol on $\partial M$ possibly up to gauge. The following is an example result of this type (compare with the scalar result in \cite[Theorem 1.4]{Oksanen2024}).

\begin{proposition} \label{prop_bd_polynomials}
Assume the conditions in Theorem \ref{thm_bd}. Suppose additionally that $m=2$, $p_m = p_m^*$, $q_m = q_m^*$, and the polynomials $z \mapsto \det(p_m(x',0,\xi',z))$ and $z \mapsto \det(q_m(x',0,\xi',z))$ have $N$ distinct roots with positive imaginary parts and that the corresponding eigenvectors are related to a simple eigenvalue $0$ and are linearly independent. Also assume that $\p M$ is noncharacteristic for $Q$ at $\hat{x}$ in the sense that $\det(q_m(\hat{x}, \nu)) \neq 0$. Then 
\begin{equation}
\det(p_m(x',0,\xi',z)) = c(x') \det(q_m(x',0,\xi',z))
\end{equation}
for $(x',\xi')$ near $(\hat{x}', \hat{\xi}')$, where $c(x')$ is a nonvanishing function.
\end{proposition}

\begin{proof}
Note that $\det(p_m(x',0,\xi',z))$ is a polynomial of degree $2N$ in $z$ with real coefficients. By assumption, there are $N$ distinct roots $z^{(1)}, \ldots, z^{(N)}$ with positive imaginary parts (the remaining roots are their complex conjugates) and corresponding linearly independent vectors $\omega^{(1)}, \ldots, \omega^{(N)}$. Similarly $\det(q_m(x',0,\xi',z))$ has distinct roots $\rho^{(s)}$ and linearly independent eigenvectors~$\beta^{(s)}$.

Varying~$\beta$, from the conclusion of Theorem \ref{thm_bd} we obtain that 
\begin{equation}
p_m(x',0,\xi',z^{(r)}) \omega^{(r)} = q_m(x',0,\xi',z^{(r)}) \omega^{(r)} = 0.
\end{equation}
But since $\det(q_m(x',0,\xi',z))$ has roots $\rho^{(s)}$ and since $\ker(q_m)$ is one-dimensional, it must follow that (after renaming indices and multiplying by scalars)  for $1 \leq r \leq N$, 
\begin{equation}
z^{(r)} = \rho^{(r)}, \qquad \omega^{(r)} = \beta^{(r)}.
\end{equation}
In particular we have, for $(x',\xi')$ near $(\hat{x}',\hat{\xi}')$ and for all $z \in \C$,  
\begin{equation}
\det(p_m(x',0,\xi',z)) = c(x',\xi') \det(q_m(x',0,\xi',z)).
\end{equation}
Note that $\lambda^{-2N} \det(p_m(x',0,\xi',\lambda)) = \det(p_m(x',0,\xi'/\lambda,1)$. Letting $\lambda \to \infty$ gives 
\begin{equation}
c(x',\xi') = c(x') = \frac{\det(p_m(x',0,0,1))}{\det(q_m(x',0,0,1))}
\end{equation}
The denominator is not zero by the noncharacteristic assumption.
\end{proof}

At this point, instead of continuing in the general setting, we specialize to the case of the elastic wave operator and give related boundary determination results. In that case we will show that the nondegeneracy conditions such as those in Proposition \ref{prop_bd_polynomials} hold generically.

\section{Boundary determination for elastic operator}
\label{sec:bdElastic}
In this section we prove a reformulated version of Lemma~\ref{lem:Alessandrini} based on the construction of special solutions in Proposition~\ref{prop_concentrating_sol}. In this section all results are formulated for the elastic wave operator~$P$
\begin{equation}
(P_{\rho,c} u)_i = \rho \p_t^2 u_i - \p_j(c_{ijkl} \p_k u_l)
\end{equation}
in $(0,T) \times \Omega$, with principal symbol $p = (\rho \xi_0^2 \delta_{il} - c_{ijkl} \xi_j \xi_k)$. We let~$x_0$ be the time variable and~$\xi_0$ the corresponding dual variable. Instead of the Cauchy data set~$C_P$, we consider boundary measurements given by the DN map 
\begin{equation}
(\Lambda_{\rho,c} f)_i = \nu_j c_{ijkl} \partial_k u_l|_{(0,T) \times \partial \Omega}
\end{equation}
defined in Section \ref{sec_intro}. Similarly we define $\tilde{P} = P_{\tilde{\rho}, \tilde{c}}$.

The DN map~$\Lambda_{\rho,c}$ involves the boundary value of the unknown tensor~$c$, unlike the Cauchy data set~$C_{P_{\rho,c}}$. Thus instead of~\eqref{bd_integral_identity} we will need to use a slightly different integral identity.

\begin{lemma} \label{dn_map_identity}
Given any $f, \tilde{f} \in C^{\infty}_c((0,T) \times \partial \Omega)$, we have 
\begin{equation}
((\Lambda_{\rho,c}-\Lambda_{\tilde{\rho}, \tilde{c}}) f, \tilde{f}) = \int_0^T \int_{\Omega} ((c_{ijkl}- \tilde{c}_{ijkl}) \partial_k u_l \overline{\partial_j \tilde{u}_i} - (\rho-\tilde{\rho}) \partial_t u_i \overline{\partial_t \tilde{u}_i}) \,dx \,dt
\end{equation}
where $P_{\rho,c} u = P_{\tilde{\rho}, \tilde{c}} \tilde{u} = 0$ with 
\begin{align}
u|_{(0,T) \times \partial \Omega} &= f, & u(0) &= \partial_t u(0) = 0, \\
\tilde{u}|_{(0,T) \times \partial \Omega} &= \tilde{f}, & \tilde{u}(T) &= \partial_t \tilde{u}(T) = 0.
\end{align}
\end{lemma}
\begin{proof}
The weak formulation of the DN map is  
\begin{equation}
(\Lambda_{\rho,c} f, \tilde{f}) = \int_0^T \int_{\Omega} (c_{ijkl} \partial_k u_l \overline{\partial_j \tilde{u}_i} - \rho \partial_t u_i \overline{\partial_t \tilde{u}_i}) \,dx \,dt
\end{equation}
where~$u$ is as in the statement of the lemma, and~$\tilde{u}$ is any smooth function with $\tilde{u}|_{(0,T) \times \partial \Omega} = \tilde{f}$ and $\tilde{u}(T) = \p_t \tilde{u}(T) = 0$. The result follows from this fact and the symmetry property $(\Lambda_{\rho,c} f, \tilde{f}) = ( f, \Lambda_{\rho,c}^T \tilde{f})$, where~$\Lambda_{\rho,c}^T$ is the analogue of $\Lambda_{\rho,c}$ but with zero Cauchy data imposed at $t=T$ instead of $t=0$. The last fact can be verified by a direct computation.
\end{proof}

Below, we will fix a boundary point on $(0,T) \times \partial \Omega$ and work in boundary normal coordinates $(x',x_n)$ so that $x' = (x_0, \ldots, x_{n-1})$ with $x_0 = t$, and the point of interest is $(\hat{x}', \hat{\xi}')$. We will also write 
\begin{equation}
\bar{c} = c - \tilde{c}, \qquad \bar{\rho} = \rho - \tilde{\rho}.
\end{equation}
Since the stiffness tensor and density are real valued, there is no confusion with complex conjugates. We will also write 
\begin{equation}
M = (0,T) \times \overline{\Omega}.
\end{equation}

\begin{proposition}\label{prop:Alessandrini}
Let $\hat{x} = (\hat{x}',0) \in (0,T) \times \partial \Omega$, and suppose that there is $\hat{\xi} = (\hat{\xi}',0) \in T_{\hat{x}}(\p M) \setminus 0$ such that the polynomials 
\begin{equation}
z \mapsto \det(p(\hat{x}', 0, \hat{\xi}', z)) \quad \text{ and } \quad \tilde{z} \mapsto \det(\tilde{p}(\hat{x}', 0, \hat{\xi}', \tilde{z}))
\end{equation}
both have $n$ distinct roots $z^{(r)}$ and $\tilde{z}^{(s)}$ with positive imaginary parts. Moreover, suppose that $\ker(p(\hat{x}', 0, \hat{\xi}', z^{(r)})) = \C \omega^{(r)}$ and $\ker(\tilde{p}(\hat{x}', 0, \hat{\xi}', \tilde{z}^{(s)})) = \C \tilde{\omega}^{(s)}$ 
. 
Suppose that the roots and eigenvectors depend smoothly on $(x',\xi')$ near $(\hat{x}', \hat{\xi}')$ as in Lemma \ref{lemma_z_omega}.

If $\Lambda_{\rho,c} = \Lambda_{\tilde{\rho},\tilde{c}}$, then for any $r, s$ we have 
\begin{equation}
(\bar{c}_{ijkl}(x',0) \zeta_k^{(r)} \overline{\tilde{\zeta}_j^{(s)}} - \bar{\rho}(x',0) \delta_{il} \zeta_0^2) \omega_l^{(r)}  \overline{\tilde{\omega}_i^{(s)}} = 0
\end{equation}
for $(x',\xi')$ near $(\hat{x}', \hat{\xi}')$. Here 
\begin{equation}
\zeta^{(r)}(x',\xi') = (\xi', z^{(r)}(x',\xi')), \qquad \tilde{\zeta}^{(s)}(x',\xi') = (\xi', \tilde{z}^{(s)}(x',\xi')).
\end{equation}
\end{proposition}
\begin{proof}
We fix roots~$z^{(r)}$ and~$\tilde{z}^{(s)}$ with corresponding eigenvectors~$\omega^{(r)}$ and~$\tilde{\omega}^{(s)}$. By Proposition \ref{prop_concentrating_sol} there are approximate solutions 
\begin{equation}
u^{(r)} = h^{1/2} e^{i\Phi^{(r)}/h} a^{(r)}, \qquad \tilde{u}^{(s)} = h^{1/2} e^{i \tilde{\Phi}^{(s)}/h} \tilde{a}^{(s)},
\end{equation}
with
\begin{equation}
\begin{split}
\nabla_{t,x} \Phi^{(r)}(x',0) &= \zeta^{(r)}
,\\
\nabla_{t,x} \tilde{\Phi}^{(s)}(x',0) &= \tilde{\zeta}^{(s)}
,\\
a^{(r)}(x',0) &= \psi(x') \omega^{(r)} + O(h)
,\\
\tilde{a}^{(s)}(x',0) &= \tilde{\psi}(x') \tilde{\omega}^{(s)} + O(h)
\end{split}
\end{equation}
with $\psi, \tilde{\psi}$ arbitrary smooth functions.
One has $\norm{u^{(r)}}_{L^2(M)} \sim h$, $\norm{u^{(r)}}_{H^1(M)} \lesssim 1$ and $\norm{Pu}_{H^1(M)} = O(h^{\infty})$ and similarly for $\tilde{u}^{(s)}$. By well-posedness we can therefore add $O_{H^1(M)}(h^{\infty})$ correction terms to~$u^{(r)}$ and~$\tilde{u}^{(s)}$ to obtain exact solutions.

As in the proof of Theorem \ref{thm_bd}, inserting such solutions in the identity of Lemma \ref{dn_map_identity} yields 
\begin{align}
0 &= \lim_{h \to 0} \int_0^T \int_{\Omega} (\bar{c}_{ijkl} \partial_k u_l \overline{\partial_j \tilde{u}_i} - \bar{\rho} \partial_t u_i \overline{\partial_t \tilde{u}_i}) \,dx \,dt \\
 &= \lim_{h \to 0} h^{-1} \int_0^T \int_{\Omega} (\bar{c}_{ijkl} \zeta^{(r)}_k \overline{\tilde{\zeta}^{(s)}_j} - \bar{\rho}\zeta^{(r)}_0 \overline{\tilde{\zeta}^{(s)}_0} \delta_{il}) \omega_l^{(r)} \overline{\tilde{\omega}^{(s)}_i} e^{i\Theta/h} \psi(x') \overline{\tilde{\psi}(x')} \,dx \,dt \\
 & \qquad + \lim_{h \to 0} h^{-1} \int_0^T \int_{\Omega} e^{i\Theta/h}(O(x_n) + O(h)) \,dx \,dt
\end{align}
where $\Theta = \Phi^{(r)} - \overline{\tilde{\Phi}^{(s)}} = (z^{(r)} - \overline{\tilde{z}^{(s)}}) x_n + O(x_n^2)$. Evaluating the limits as in Theorem \ref{thm_bd} and varying~$\psi$ and~$\tilde{\psi}$ gives that 
\begin{equation}
(\bar{c}_{ijkl}(x',0) \zeta_k^{(r)} \overline{\tilde{\zeta}_j^{(s)}} - \bar{\rho}(x',0) \delta_{il} \zeta_0^2) \omega_l^{(r)}  \overline{\tilde{\omega}_i^{(s)}} = 0
\end{equation}
for $(x',\xi')$ near $(\hat{x}', \hat{\xi}')$.
\end{proof}

\section{Algebraic results}
\label{sec:algebra}

In this section we prove the algebraic results in Lemma~\ref{lem:genPair}, Remark~\ref{rem:failureset} and give the full proof of Theorem~\ref{thm:main} that is here formulated as Theorem~\ref{thm:generic-uniqueness}. For this, we need some notions and results from algebraic geometry that we will introduce now. An invitation to algebraic geometry specifically for analysts is offered in~\cite{Smith2000Invitation}. For a more thorough algebraic introduction, we refer to~\cite{atiyah1969introduction}. The algebraic geometry background can be found in~\cite{Hartshorne}.

\subsection{Notion and results from algebraic geometry.}
 The ambient spaces of our geometric objects will be product spaces of affine space~$\C^n$ and projective space $\CP^{n-1} = (\C^n \setminus \{0\})/\sim$, where $\omega_1 \sim \omega_2$ if and only if $\omega_1 = \lambda \omega_2$ for some~$\lambda \in \C\setminus \{0\}$. By $\C[x_1,\ldots,x_n]$ and $\C(x_1,\ldots,x_n)$, we denote the ring of polynomials and the field of rational functions, respectively, in~$n$ indeterminates with complex coefficients. For a \textit{homogeneous} polynomial $f \in \C[x_1,\ldots,x_n]$ and $\omega \in \CP^{n-1}$, the relation $f(\omega)=0$ is well-defined, that is, independent of the choice of representative from~$\C^n$ of~$\omega$. In the following, we will introduce notions only for~$\C^n$, and the corresponding notions for~$\CP^{n-1}$ can be inferred by using homogeneous polynomials when necessary.

A non-empty subset $I \subseteq \C[x_1,\ldots,x_n]$ such that $f-g \in I$ and $rf \in I$, for all $f,g \in I$ and $r \in \C[x_1,\ldots,x_n]$, is called an \textit{ideal}.
A subset $S \subseteq \C[x_1,\ldots,x_n]$ of polynomials defines the \textit{algebraic set} $V(S) \subseteq \C^n$ of common zeros of elements in~$S$. It can be easily seen that $V(S) = V(\langle S \rangle)$, where $\langle S \rangle$ is the smallest ideal containing~$S$. \textit{Hilbert's Basis Theorem} says that, for every ideal $I \subseteq \C[x_1,\ldots,x_n]$, there exists a finite set $S = \{f_1,\ldots,f_m\}$ such that $I = \langle S \rangle = \langle f_1,\ldots,f_m\rangle$. In particular, the algebraic sets are exactly the solution sets to \textit{finite} systems of polynomial equations. Moreover, the algebraic sets are the closed sets of a topology on~$\C^n$, the \textit{Zariski topology}, which is coarser than the Euclidean topology.\\ 

\noindent
\textbf{Algebro-geometric correspondences.} An ideal~$I$ is \textit{radical} if $f^s \in I$ implies $f \in I$ for every $f \in \C[x_1,\ldots,x_n]$ and every positive integer~$s$. An ideal~$P$ is \textit{prime} if $fg \in P$ implies $f \in P$ or $g \in P$ for all $f,g \in \C[x_1,\ldots,x_n]$. Every \textit{maximal} ideal, that is, every ideal maximal among ideals that are strictly contained in $\C[x_1,\ldots,x_n]$, is prime. Every prime ideal is radical. We can define the analogous types of ideals for any ring~$R$. The set of prime ideals of~$R$ is denoted by~$\Spec R$.

A subset of a topological space is \textit{irreducible} if it is not the union of two closed subsets properly contained in it. \textit{Hilbert's Nullstellensatz} says that the map $I \mapsto V(I)$ gives a bijective correspondence 
\begin{equation}
\text{radical ideals} \longleftrightarrow \text{algebraic sets}
\end{equation}
whose inverse sends a subset of~$\C^n$ to the set of all polynomials vanishing on it and who restricts to bijective correspondences
\begin{equation}
\text{prime ideals} \longleftrightarrow \text{irreducible algebraic sets}
\end{equation}
and
\begin{equation}
\text{maximal ideals} \longleftrightarrow \text{points in } \C^n.
\end{equation}
Using the existence of primary decompositions in polynomial rings, one can prove that every algebraic set is decomposed uniquely as a finite union of its maximal irreducible subsets, also called \emph{irreducible components}.\\

\noindent
\textbf{Rings of functions and morphisms.} Let $V \subseteq \C^n$ be an algebraic set with the Zariski topology, that is, the subspace topology coming from the Zariski topology on~$\C^n$. We endow~$V$ with extra algebraic structure, namely, with its \textit{sheaf of regular functions}: For every open subset $U \subseteq V$, we denote by~$\mathcal{O}_V(U)$ the ring of rational functions on~$U$. The topological space~$V$ together with this sheaf is called an \textit{affine algebraic variety} (over~$\C$). In a similar way, one can endow algebraic sets in~$\CP^{n-1}$ with a sheaf structure, which leads to the concept of a \textit{projective algebraic variety}.

Morphisms between algebraic varieties are continuous maps that preserve the sheaf structure. As~$\CP^{n-1}$ can be covered by open subsets that are (as algebraic varieties) isomorphic to the algebraic variety~$\C^n$, every algebraic variety can be covered by open subsets isomorphic to affine algebraic varieties.
Furthermore, we have an equivalence of categories
\begin{align}\label{equ:equivalence}
    \text{affine algebraic varieties } \longleftrightarrow& \text{ finitely generated $\C$-algebras}
\end{align}
given by $V \mapsto \mathcal{O}_V(V)$, where all our algebras are commutative and associative with a unit element. \\

\noindent
\textbf{Schemes.} To arrive at the language of modern algebraic geometry that we need to use in some of our arguments, we have to enrich an algebraic variety~$V$ in one further step. By Hilbert's Nullstellensatz and the equivalence from (\ref{equ:equivalence}), we can identify the points in~$V$ with the maximal ideals of~$\mathcal{O}_V(V)$. Now, instead of just taking the maximal ideals, we can use the space~$\Spec \mathcal{O}_V(V)$ of \textit{all} prime ideals. The points of the variety~$V$ now correspond to the \textit{closed points} of~$\Spec \mathcal{O}_V(V)$.
The Zariski topology and the sheaf structure extend naturally and we call the resulting object the \textit{affine scheme} associated to~$V$. 

There is a canonical way to define a Zariski topology and a sheaf structure on~$\Spec R$ for any commutative ring~$R$, which leads us to the general notion of an \textit{affine scheme}. A \textit{scheme} is now a topological space with a (suitable) sheaf structure that can be covered by open subsets that are isomorphic to affine schemes.
By gluing affine open sets, every algebraic variety (over~$\C$) has an associated scheme (together with a morphism of schemes to~$\Spec \C$), which gives us a fully faithful functor
\begin{align}\label{equ:scheme-functor}
    \text{Varieties over } \C \longrightarrow \text{Schemes with a morphism to } \Spec \C.
\end{align} 
Hence, we can identify a variety with its associated scheme. For instance, by abuse of notation, we denote by~$\C^n$ and~$\CP^{n-1}$ the schemes associated to these spaces, respectively.

As an alternative definition (that is equivalent up to isomorphisms), one can now say that an algebraic variety (over~$\C$) is a \textit{reduced} scheme~$V$ together with a morphism of schemes $V \to \Spec \C$  that is \textit{separated} and \textit{of finite type}. Along with others, these technical notions will be used in the proof of Theorem~\ref{thm:generic-uniqueness} and are explained in the following paragraph. \\

\noindent
\textbf{Properties of (morphisms of) schemes.} A ring~$R$ is \textit{reduced} if it has no non-zero nilpotent elements. $R$ is an \textit{integral domain} if it has no non-zero zero divisors. A scheme~$V$ is \textit{reduced} or \textit{integral}, respectively, if it can be covered by open subsets isomorphic to affine schemes~$\Spec R$, where~$R$ varies over a set of reduced rings or integral domains, respectively.

A ring homomorphism $\varphi: A \to B$ induces a morphism of affine schemes $\Spec B \to \Spec A$ by $P \mapsto \varphi^{-1}(P)$.
A morphism of schemes $f: X \to S$ is of \textit{finite type} if it is locally of the form $\Spec B \to \Spec A$, where~$B$ is a finitely generated $A$-algebra via a map $A \to B$.

The way in which a morphism is separated is described most intuitively in the following way, though this is no formal definition.
Morphisms $X \to S$, where~$S$ is a fixed scheme, form the category of $S$-schemes. A morphism $f: X \to S$ is \textit{separated} if it satisfies a condition in the category of $S$-schemes that is analogous to the Hausdorff property in the category of topological spaces. Note that the underlying topological space of a scheme that consists of more than one point is never Hausdorff. 

Let $X\to S$ and~$Y \to S$ be $S$-schemes. The \textit{fiber product} $X \times_S Y$ is their pullback in the category of schemes. For instance, if $X \to \Spec \C$ and $Y \to \Spec \C$  are varieties over $\C$ then $X \times_{\Spec \C} Y$ is just (the scheme associated to) the usual product of varieties. Moreover, if $X = \Spec A$, $Y = \Spec B$ and $S = \Spec R$ are affine schemes then $X \times_S Y$ is isomorphic to the affine scheme associated to $A \otimes_R B$.

A morphism $f:X \to S$ is \textit{flat} if it is locally of the form $\Spec B \to \Spec A$, where $B$ is a flat $A$-algebra, that is, tensoring $A$-linear maps of $A$-modules with~$B$ preserves injectivity.

A morphism $f: X \to S$ is \emph{smooth} at a point $x \in X$ if, around $x$, it is locally of the form $\Spec A[x_1,\ldots,x_n]/\langle f_1,\ldots,f_m\rangle \Spec A$ such that the Jacobian matrix $(\partial f_i/\partial x_j)$ is of rank~$m$ at~$x$.\footnote{The map of spectra is just the one induced by the natural ring map $A \to A[x_1,\ldots,x_n]/\langle f_1,\ldots,f_m\rangle$ via (\ref{equ:equivalence}) and (\ref{equ:scheme-functor}), where the latter is the ring of equivalence classes of polynomials over~$A$ subject to the relations $f \sim g$ if and only if $f-g \in \langle f_1,\ldots,f_m \rangle$.} 

Note, that a morphism is always smooth in the sense of manifold theory because it is locally defined by algebraic functions. However, not every morphism is a smooth morphism in the sense of algebraic geometry. For instance, if~$X$, $S$ and~$f$ are defined over~$\C$ and the sets of complex points~$X(\C)$ and~$S(\C)$ form smooth complex manifolds, then~$f$ is a smooth morphism if and only if the differential~$Df$ of the induced map $f: X(\C) \to S(\C)$ of smooth manifolds is surjective. Moreover, every smooth morphism of schemes is flat.

The \emph{generic point}~$\eta$ of an integral scheme~$S$ lies in every affine open cover set $\Spec A \subset S$ and $\eta = \langle 0 \rangle \subset A$ is the zero ideal. If $X \to S$ is a morphism then the \emph{generic fiber} $X_\eta$ is just the fiber over the generic point. One main feature of our proofs, which we will explicitly use in the proofs of Lemma~\ref{lem:squareFree} and Lemma~\ref{lem:generic-noncontainment}, is that some properties of the generic fiber~$X_\eta$ spread out to fibers~$X_s$ for~$s$ in a non-empty Zariski open and hence dense subset of~$S$. That is, some properties of the generic fiber are generic properties of the fibers.






\subsection{Proof of the main result Theorem~\ref{thm:main}.}

In this section we will use the following notation that was partly introduced in Section~\ref{sec:proofidea}.
Denote by
\begin{equation}
\Xi
=
\{
(\xi,\tilde\xi)
\in(\CP^n)^2
;
\xi_i=\tilde\xi_i
\text{ for all }
i=0,\dots,n-1
\}
\end{equation}
the space of spatial covectors.
The space of pairs of models will be
\begin{equation}
M_\K
=
(\K\times\K^N)^2
\end{equation}
with $\K\in\{\R,\C\}$.
We consider a pair of models
\begin{equation}
m
=
(\rho,c,\tilde\rho,\tilde c)
\in
M_\C
\end{equation}
to be a parameter.
Given such a parameter,
we define
in the ambient space
\begin{equation}
A
=
\Xi\times(\CP^n)^2
\end{equation}
the polynomials
\begin{equation}
\begin{split}
F(\xi,\tilde\xi,\omega,\tilde\omega)
&=
[
(c-\tilde c)_{ijkl}
\xi_k\tilde\xi_j
-
(\rho-\tilde\rho)
\delta_{il}
\xi_0\tilde\xi_0
]
\omega_l
\tilde\omega_i
,\\
f_i(\xi,\tilde\xi,\omega,\tilde\omega)
&=
(
c_{ijkl}\xi_j\xi_k
-
\rho\delta_{il}\xi_0^2
)
\omega_l
,\\
\tilde f_i(\xi,\tilde\xi,\omega,\tilde\omega)
&=
(
\tilde c_{ijkl}\tilde\xi_j\tilde\xi_k
-\tilde\rho\delta_{il}\tilde\xi_0^2
)
\tilde\omega_l
\end{split}
\end{equation}
for $i=1,\dots,n$.
These depend on $m$ but we leave it out of the notation.
Let us denote the principal symbols by $p = (
c_{ijkl}\xi_j\xi_k
-
\rho\delta_{il}\xi_0^2
)_{il}$ and $\tilde p = (
\tilde c_{ijkl}\tilde\xi_j\tilde\xi_k
-\tilde\rho\delta_{il}\tilde\xi_0^2
)_{il}$.

With these, we define two varieties in~$A$:
\begin{equation}
\begin{split}
W_m
&=
\{
(\xi,\tilde\xi,\omega,\tilde\omega)
\in A
;
\\\quad
f_i&=0
\text{ and }
\tilde f_i=0
\text{ for all }
i=1,\dots,n
\}
\end{split}
\end{equation}
and
\begin{equation}
U_m
=
\{
(\xi,\tilde\xi,\omega,\tilde\omega)
\in A
;
F
=
0
\}
.
\end{equation}
The variety~$W_m$ simply describes points on the complex slowness surfaces and the associated polarizations of the two models when the cotangent vectors only differ in the last component.
The variety~$U_m$ describes the additional constraint that we obtain through the Alessandrini identity and special solutions.

We will also make use of the realish subset
\begin{equation}
A_\R
=
\{
(\xi,\tilde\xi,\omega,\tilde\omega)
\in A
;
\xi_i\in\R
\text{ for all }
i=0,\dots,n-1
\}
.
\end{equation}

Note that the polynomials $F,f_i,\tilde f_i$ are homogeneous and, hence, define a variety in the multiprojective space~$A$.

\begin{definition}
We say that $(\rho,c)\in\C\times\C^N$ is \emph{squarefree} if $\xi\mapsto\det(p(\xi))$ is squarefree in $\C[\xi_0,\dots,\xi_n]$.
We say that a model $m=(\rho,\tilde\rho,c,\tilde c)\in M_\C$ is squarefree if both $(\rho,c)$ and $(\tilde\rho,\tilde c)$ are squarefree.
\end{definition}

In the following Lemma we show that the slowness polynomial $\det(p(\xi))$ being squarefree in $\C[\xi_0,\dots,\xi_n]$ is equivalent to the $\det(p(\xi))$ being squarefree in $\C[\xi_n]$ and thus~$\det(p(\xi))$ having simple roots in~$\xi_n$:

\begin{lemma}\label{lem:squareFree}
    Let~$\sigma$ be a square-free polynomial in $\C[\xi_0,\ldots,\xi_{n-1},\xi_n]$. Then 
    \begin{equation}
    \{(a_0,\ldots,a_{n-1}) \in \C^n \mid \sigma(a_0,\ldots,a_{n-1},\xi_n) \text{ is square-free in } \C[\xi_n]\}
    \end{equation}
    contains a non-empty Zariski open subset of~$\C^n$, for arbitrary dimensions~$n$.
\end{lemma}
\begin{proof}[Proof of Lemma \ref{lem:squareFree}.]
    We use that a polynomial is (by definition) square-free if and only if the affine scheme  defined by it is reduced.
    The two unique factorization domains $\C[\xi_0,\ldots,\xi_{n-1},\xi_n]$ and $\C[\xi_0,\ldots,\xi_{n-1}][\xi_n]$ are (canonically) isomorphic, and any ring isomorphism preserves irreducibility of elements. Therefore, we may consider~$\sigma$ as a squarefree polynomial in~$\xi_n$ with coefficients in $\C[\xi_0,\ldots,\xi_{n-1}]$, that is, $\sigma = g_1 \cdots g_k$, where the $g_i \in \C[\xi_0,\ldots,\xi_{n-1}][\xi_n]$ are irreducible polynomials that are pairwise non-associated (that is, $g_i \neq kg_j$ for all $k \in \C$ if $i \neq j$). The straightforward direction of Gauss' Lemma shows that $g_i \in \C(\xi_0,\ldots,\xi_{n-1})[\xi_n]$ is either invertible or irreducible. So $\sigma$ is square-free in $\C(\xi_0,\ldots,\xi_{n-1})[\xi_n]$.

    As $K = \C(\xi_0,\ldots,\xi_{n-1})$ has characteristic~$0$, this is equivalent to saying that $\gcd_K(p,p') = 1$, where $\sigma'$ is the $\xi_n$-derivative of $\sigma$. Equivalently, the ideal $(\sigma,\sigma')_{K[\xi_n]}$ equals $K[\xi_n]$. As $K$ is a field, its algebraic closure $\overline{K}$ is a $K$-vector space, that is, a free (and, hence, flat) $K$-module. So, by the very definition of flatness, tensoring the equality $(\sigma,\sigma')_{K[\xi_n]} = K[\xi_n]$  with~$\overline{K}$ over~$K$ gives an equality:
    \begin{equation}
    (\sigma,\sigma')_{\overline{K}[\xi_n]} = (\sigma,\sigma')_{K[\xi_n]} \otimes_K \overline{K} = K[\xi_n] \otimes_K \overline{K} = \overline{K}[\xi_n].
    \end{equation}
    So $\gcd_{\overline{K}}(\sigma,\sigma') = 1$ and~$\sigma$ is square-free in~$\overline{K}[\xi_n]$.

    Consider the variety
    \begin{equation}
    V = \{ (a_1,\ldots,a_n,c) \in \A^{n+1}_\C \mid p(a_1,\ldots,a_n,c) = 0 \}
    \end{equation}
    and the morphism $f = \pi|_{V}: V \to \A^n_\C$, where~$\pi$ is just the projection to the coordinates $a_1,\ldots,a_n$ of $\A^{n+1}_\C$. The fiber~$V_\eta$ of the generic point~$\eta$ of~$\A^n_\C$ has underlying variety
    \begin{equation}
    \{ z \in \A^1_{K} \mid p(\xi_0,\ldots,\xi_{n-1},\xi_n) = 0\}
    \end{equation}
    What we just showed, namely that~$\sigma$ is square-free in~$\overline{K}[\xi_n]$, means nothing else than~$V_\eta$ being geometrically reduced, that is, the scheme defined by~$\sigma$ in $\A^1_{\overline{K}}$ being reduced. Clearly,~$f$ is of finite type, so \cite[Lemma 37.26.4]{stacks-project} shows the existence of the desired non-empty Zariski open set.
\end{proof}

Conversely, simple roots for a single choice of the remaining variables already give squarefreeness:

\begin{lemma}[Simple roots at one point]\label{lem:simple-roots-sqf}
Let $\sigma\in\C[\xi_0,\dots,\xi_n]$ and $j\in\{0,\dots,n\}$, and suppose that the coefficient of the highest power of~$\xi_j$ in~$\sigma$ is a non-zero constant. For $a=(a_k)_{k\neq j}\in\C^n$, denote by $\sigma(a,\xi_j)$ the polynomial in~$\xi_j$ obtained by substituting $\xi_k=a_k$ for all $k\neq j$. If $\sigma(a,\xi_j)$ has only simple roots in $\xi_j$ for some~$a$, then~$\sigma$ is square-free in $\C[\xi_0,\dots,\xi_n]$.
\end{lemma}

\begin{proof}
Suppose that $\sigma=g^2h$ with a non-constant $g\in\C[\xi_0,\dots,\xi_n]$. The leading coefficient of~$\sigma$ in~$\xi_j$ is the product of the leading coefficients of~$g^2$ and~$h$ in~$\xi_j$. Since it is a non-zero constant, the leading coefficient of~$g$ in~$\xi_j$ is a non-zero constant as well. Hence~$g$ has positive degree in~$\xi_j$, as otherwise~$g$ would be constant, and this degree does not drop when $\xi_k$ is replaced by~$a_k$ for $k\neq j$. Thus $g(a,\xi_j)$ has a root~$r$, and since $\sigma(a,\xi_j)=g(a,\xi_j)^2h(a,\xi_j)$, the number~$r$ is a root of multiplicity at least two of $\sigma(a,\xi_j)$, a contradiction.
\end{proof}

To ensure that the set of models $(\rho,c)$ for which $\det (p(\xi',z))$ is squarefree in~$z$ contains a Zariski open subset that is non-empty (see Lemma~\ref{lem:generic-simplicity}(1)) we construct the following examples for squarefreeness in general dimensions:
\begin{lemma}[Examples for squarefreeness]\label{lem:squarefreeness}
For $n\geq2$ let $c^*_{ijk\ell}=\min\{i,j\}\,(\delta_{ik}\delta_{j\ell}+\delta_{i\ell}\delta_{jk})$ and $\rho_*=1$. Such a stiffness tensor reads as follows in dimensions $n=2,3$ in Voigt notation:
\begin{align}
c^*&=\begin{pmatrix}
    c_{1111} & c_{1122} & c_{1112}\\
    & c_{2222} & c_{1222}\\
    \mathclap{\text{s.}} &  & c_{1212}
\end{pmatrix}=\begin{pmatrix}2&0&0\\0&4&0\\0&0&1\end{pmatrix}
\quad(n=2),\\
c^*&=\begin{pmatrix} c_{1111} & c_{1122} & c_{1133} & c_{1123} & c_{1113} & c_{1112}\\
    & c_{2222} & c_{2233} & c_{2223} & c_{2213} & c_{2212}\\
    & & c_{3333} & c_{3323} & c_{3313} & c_{3312}\\
    &  &  & c_{2323} & c_{2313} & c_{2312}\\
    &  \mathclap{\text{symm.}} &  &  & c_{1313} & c_{1312}\\
    &  &  &  &  & c_{1212}
    \end{pmatrix} =\begin{pmatrix}
2&0&0&0&0&0\\
0&4&0&0&0&0\\
0&0&6&0&0&0\\
0&0&0&2&0&0\\
0&0&0&0&1&0\\
0&0&0&0&0&1
\end{pmatrix}
\quad(n=3).
\end{align}
Then $c^*$ satisfies \eqref{eq:sym} and \eqref{eq:pos}, the model $(\rho_*,c^*)$ is squarefree, and for $\xi'_*=(1,0,\dots,0)$ all roots of $z\mapsto\det p(\xi'_*,z)$ are simple.
\end{lemma}

\begin{proof}
The tensor $c^*$ satisfies \eqref{eq:sym}, and its Voigt matrix is diagonal with positive entries $c^*_{iiii}=2i$ and $c^*_{ijij}=\min\{i,j\}$ for $i\neq j$. The positivity condition in \eqref{eq:pos} is equivalent to its Voigt matrix being positive definite as can be seen from the following: Enumerate index pairs $(i,j)$ with $i\leq j$ by $I=1,\dots,n(n+1)/2$, and express a symmetric $n\times n$ matrix~$A$ as the vector $\hat A\in\R^{n(n+1)/2}$ that satisfies $\hat A_I=A_{ii}$ if $I$ corresponds to $(i,i)$ and $\hat A_I=2A_{ij}$ if $I$ corresponds to $(i,j)$ with $i<j$. Then
\begin{equation}
        c_{ijk\ell}A_{ij}A_{k\ell}
        =
        \hat A^Tc \hat A,
\end{equation}
where $c$ on the right-hand side is expressed in Voigt notation. Note that by~\eqref{eq:sym} the terms with $i\neq j$ or $k\neq\ell$ occur twice in the sum on the left. Hence, as the Voigt matrix is diagonal with positive entries $c^*$ is positive definite. Since $\Gamma_{c^*}(\xi)_{i\ell}=\min\{i,\ell\}\,\xi_i\xi_\ell+\delta_{i\ell}\sum_{j}\min\{i,j\}\,\xi_j^2$, we have $\Gamma_{c^*}(e_n)=\operatorname{diag}(1,2,\dots,n-1,2n)$ and $p(\xi'_*,z)=z^2\,\Gamma_{c^*}(e_n)-I$. Its determinant $(2nz^2-1)\prod_{i=1}^{n-1}(iz^2-1)$ has the $2n$ distinct roots $\pm1,\pm2^{-1/2},\dots,\pm(n-1)^{-1/2},\pm(2n)^{-1/2}$. Every entry of $p(\xi)$ has degree at most two in~$\xi_n$, and the coefficient of~$\xi_n^2$ in $p(\xi)$ is $\Gamma_{c^*}(e_n)$. Hence the coefficient of the highest power $\xi_n^{2n}$ in $\det p(\xi)$ is the non-zero constant $\det\Gamma_{c^*}(e_n)=2n\,(n-1)!$, and Lemma~\ref{lem:simple-roots-sqf} with $j=n$ gives that $(\rho_*,c^*)$ is squarefree.
\end{proof}

\begin{definition}
We say that $m=(\rho,\tilde\rho,c,\tilde c)\in M_\R$ is a \emph{valid pair} if it is squarefree and
there is a non-empty open set $O\subset A_\R\cap W_m$ in the Euclidean topology with $O\subset U_m$ that contains a point at which $\xi_n$ and $\tilde\xi_n$ are simple roots of $\det p(\xi_0,\dots,\xi_{n-1},\,\cdot\,)$ and $\det\tilde p(\tilde\xi_0,\dots,\tilde\xi_{n-1},\,\cdot\,)$, respectively.
\end{definition}

The key result from analysis is this:

\begin{lemma}[Algebraic interpretation of Alessandrini]\label{lem:alg-alessandrini}
If $\Lambda_{\rho,c}=\Lambda_{\tilde \rho, \tilde c}$ and for some $x\in\partial\Omega$ the pair $(\rho(x),\tilde\rho(x),c(x),\tilde c(x))$ of models is squarefree, then it is a valid pair.
\end{lemma}

The proof of Lemma~\ref{lem:alg-alessandrini} is given after Lemma~\ref{lem:generic-simplicity}.

\begin{lemma}[Generic non-containment]\label{lem:generic-noncontainment}
The set
\begin{equation}
G=\{m\in M_\C:\text{no irreducible component of $W_m$ is contained in $U_m$}\}
\end{equation}
contains a non-empty Zariski open subset of~$M_\C$. 
\end{lemma}

\begin{proof}
Write $\zeta=(\xi_1,\ldots,\xi_n)$ and
$\tilde\zeta=(\tilde\xi_1,\ldots,\tilde\xi_n)$, and let
\begin{equation}
M^\circ=\left\{m\in M_\C:\rho\tilde\rho
\det(c_{inn\ell})_{i,\ell=1}^n
\det(\tilde c_{inn\ell})_{i,\ell=1}^n\ne 0\right\}.
\end{equation}
This is a non-empty Zariski open subset of~$M_\C$. We first show that the
total family of all the~$W_m$ over~$M^\circ$ is integral.

Consider the open subset $B\subset A$ defined by
$\xi'=(\xi_0,\ldots,\xi_{n-1})\ne0$ and
$\zeta\ne0\ne\tilde\zeta$. On this set the common homogeneous coordinates
of~$\xi$ and~$\tilde\xi$ can be written as
\begin{equation}
[\xi_0:\cdots:\xi_{n-1}:\xi_n:\tilde\xi_n]\in\CP^{n+1}.
\end{equation}
Thus~$B$ is an open subset of
$\CP^{n+1}\times\CP^{n-1}\times\CP^{n-1}$, and is irreducible.
For $m\in M^\circ$, the equations defining~$W_m$ have no solutions outside
$B$. Indeed, if $\zeta=0$, then the first system reduces to
$-\rho\xi_0^2\omega=0$, contradicting $\rho\ne0$. The same argument
applies if $\tilde\zeta=0$. If $\xi'=0$, both~$\xi$ and~$\tilde \xi$ are
multiples of~$e_n$, and the invertibility of
$(c_{inn\ell})$ and $(\tilde c_{inn\ell})$ excludes non-zero kernel
vectors.

Let $\mathcal W\subset M_\C\times B$ be the scheme defined by
$f_i=\tilde f_i=0$ for $i=1,\ldots,n$. For fixed
$(\xi,\tilde\xi,\omega,\tilde\omega)\in B$, these equations are linear in
the coefficients of the two models. We claim that they impose~$2n$
independent linear conditions. In fact, varying~$c$ among the symmetric
rank-one tensors $c_{ijk\ell}=v_iv_jv_kv_\ell$ gives the $c$-dependent part of~$f_i$ as
\begin{align}
c_{ijk\ell}\zeta_j\zeta_k\omega_\ell
&=v_iv_jv_kv_\ell\zeta_j\zeta_k\omega_\ell\\
&=v_i\left(\sum_j v_j\zeta_j\right)
       \left(\sum_k v_k\zeta_k\right)
       \left(\sum_\ell v_\ell\omega_\ell\right)\\
&=v_i(v\cdot\zeta)^2(v\cdot\omega).
\end{align}
in the first system. These span~$\C^n$, as otherwise a non-zero vector~$a$
would satisfy
\begin{equation}
(a\cdot v)(v\cdot\zeta)^2(v\cdot\omega)=0
\quad\text{for all }v\in\C^n,
\end{equation}
which is impossible since~$a$,~$\zeta$ and~$\omega$ are non-zero and hence the polynomials $v\cdot a$, $v\cdot \zeta$ and $v\cdot \omega$ in~$v$ are non-zero. The same
argument applies to~$\tilde c$, and the two sets of coefficients vary
independently. Hence $\mathcal W\to B$ is the kernel of a surjective
morphism of vector bundles and is, therefore, itself a vector bundle over the integral scheme~$B$. In particular,
$\mathcal W$ is integral.

Put $\mathcal W^\circ=\mathcal W\cap(M^\circ\times B)$. To see that
$\mathcal W^\circ\to M^\circ$ is surjective, fix $m\in M^\circ$ and consider
\begin{equation}
\det\big(c_{ijk\ell}\xi_j\xi_k-\rho\delta_{i\ell}\xi_0^2\big)=0,
\qquad
\det\big(\tilde c_{ijk\ell}\tilde\xi_j\tilde\xi_k
             -\tilde\rho\delta_{i\ell}\xi_0^2\big)=0.
\end{equation}
By the choice of~$M^\circ$, these are two non-zero homogeneous polynomials of degree~$2n$ in
$\xi_0,\ldots,\xi_{n-1},\xi_n,\tilde\xi_n$. They have a common
projective zero, since $n\geq2$, see
\cite[Theorem I.7.2]{Hartshorne}. By the choice of~$M^\circ$, this zero lies
in the part where~$\xi'$, $\zeta$ and~$\tilde\zeta$ are non-zero. Choosing
non-zero vectors in the two kernels gives a point of~$W_m$.
Consequently,~$\mathcal W^\circ$ is a non-empty open subset of
the integral scheme~$\mathcal W$, and the projection
$\mathcal W^\circ\to M^\circ$ is dominant, that is, its image is Zariski dense in the target space. Since the generic
fiber of a dominant morphism of integral schemes is integral,
the generic fiber of $\mathcal W^\circ\to M^\circ$ is non-empty
and integral.

It remains to show that~$F$ does not vanish identically on
$\mathcal W^\circ$. For this, choose the same isotropic tensor for both
models,
\begin{equation}
c_{ijk\ell}=\tilde c_{ijk\ell}
=\delta_{ik}\delta_{j\ell}+\delta_{i\ell}\delta_{jk},
\qquad \rho=2,\quad\tilde\rho=3.
\end{equation}
Let $e_1,e_n\in\C^n$ be the first and last standard basis vectors, and set
\begin{equation}
\xi_0=\tilde\xi_0=1,\qquad
\zeta=\omega=e_1,\qquad
\tilde\zeta=\tilde\omega=e_1+\frac{1}{\sqrt{2}}e_n.
\end{equation}
The two covectors agree in the components $0,\ldots,n-1$, as required.
For this tensor, the matrix $c_{ijk\ell}\zeta_j\zeta_k$ is
$\zeta\zeta^T+(\zeta\cdot\zeta)I$. Thus
$f_i=\tilde f_i=0$ for all~$i$, while
\begin{equation}
F=-(\rho-\tilde\rho)\xi_0^2\omega\cdot\tilde\omega=1.
\end{equation}

Moreover, the matrix $c_{inn\ell}$ with entries $(c_{inn\ell})_{i\ell}=(\tilde c_{inn\ell})_{i\ell}=(I+e_ne_n^T)_{i\ell}$ is invertible,
so this point belongs to~$\mathcal W^\circ$. It follows that~$F$ is not
identically zero on the integral generic fiber of
$\mathcal W^\circ\to M^\circ$.

The zero set of~$F$ is therefore nowhere dense in the generic fiber. By
\cite[Lemma 37.24.3]{stacks-project}, after replacing~$M^\circ$ by a
non-empty open subset~$K$, the intersection $W_m\cap U_m$ is nowhere
dense in~$W_m$ for all $m\in K$. Hence no irreducible component of
$W_m$ is contained in~$U_m$ for these parameters.
\end{proof}

For a model $(\rho,c)$, we write 
\begin{equation}
p(\xi)_{i\ell} = c_{ijkl}\xi_j\xi_k - \rho\delta_{il}\xi_0^2.
\end{equation}
In other words, $p_{i\ell}$ is the coefficient of~$f_i$ at~$\omega_\ell$ specialized at $(\rho,c)$.

\begin{lemma}[Generic simplicity of roots]\label{lem:generic-simplicity}
    Suppose that there exists a model $(\rho_*,c_*)$ and $\xi'_* = (\xi_{0,*}, \ldots, \xi_{n-1,*})$ such that all roots of $z \mapsto p_*(\xi'_*,z)$ are simple, whose existence is secured by Lemmas~\ref{lem:squarefreeness}, \ref{lem:simple-roots-sqf} and~\ref{lem:squareFree}. Then: 
    \begin{enumerate}
        \item The set of models $(\rho,c)$ for which there exists~$\xi'$ such that all roots of $z \mapsto \det p(\xi',z)$ are simple contains a non-empty Zariski open subset of the space of models. 
        \item There exists a non-empty Zariski open subset $M^\mathrm{simp} \subset M_\C$ such that, for every $ m = ((\rho,c),(\tilde\rho, \tilde c)) \in M^{\mathrm{simp}}$, the set
        \begin{equation}
        \Xi_m^\mathrm{simp} = \{ \xi' \in \CP^{n-1} : \text{ all roots of } z \mapsto \det p(\xi',z) \text{ and } z \mapsto \det \tilde p(\xi',z) \text{ are simple}\}
        \end{equation}
        is non-empty Zariski open in~$\CP^{n-1}$.
    \end{enumerate}
\end{lemma}

\begin{proof}
    To prove (1), we first fix a model $(\rho,c)$. The polynomial $\det p(\xi',z)$ has degree at most~$2n$ and its coefficient at~$z^{2n}$ is $a(c) = \det(c_{inn\ell})_{i,\ell = 1}^n$. Therefore, independently of~$\xi'$, whenever $a(c) \neq 0$, this polynomial has degree exactly~$2n$. 

    Let $D(\rho,c,\xi') = \disc_z \det p(\xi',z)$ be its discriminant with respect to~$z$, that is, the defining polynomial of the Zariski closure of the set of critical values of the projection from $(\xi',z)$-space to $z$-space. Consequently, if $a(c) \neq 0$, then $D(\rho,c,\xi') \neq 0$ if and only if all~$2n$ zeros of $\det p(\xi',z)$ are simple. 

    By assumption, there exists $(\rho_*,c_*,\xi'_*)$ such that $D(\rho_*,c_*,\xi'_*) \neq 0$ and, hence, the polynomial $(\rho,c) \mapsto D(\rho,c,\xi'_*)$ is not identically zero. It follows that 
    \begin{equation}
    V_* = \{ (\rho,c) : \rho a(c) D(\rho,c,\xi'_*)\neq 0 \}
    \end{equation}
    Is a non-empty Zariski open subset of the space of models. For every $(\rho,c) \in V_*$, all roots of $\det p(\xi'_*,z)$ are simple, which proves the first assertion.

    For (2), fix $(\rho,c) \in V_*$. Since $D(\rho,c,\xi'_*) \neq 0$, the polynomial $\xi' \mapsto D(\rho,c,\xi')$ is not identically zero. Hence, its non-vanishing locus is non-empty and Zariski open in~$\CP^{n-1}$. (Note for this that eliminating variables from a homogeneous ideal yields a homogeneous ideal and, hence, $\det p(\xi',z)$ is homogeneous.) 

    We now go to pairs of models. Define $M^\mathrm{simp} = V_* \times V_*$ which non-empty and Zariski open in~$M_\C$. For any $m = ((\rho,c),(\tilde \rho,\tilde c)) \in M^\mathrm{simp}$ both discriminants are non-zero at our fixed~$\xi'_*$. Therefore, the polynomial $\xi' \mapsto D(\rho,c,\xi')D(\tilde \rho,\tilde c,\xi')$ is not identically zero. It follows that the set $\Xi^\mathrm{simp}_m$ from the lemma is non-empty and Zariski open.
\end{proof}

Part~(1) above also follows from Lemma~\ref{lem:squareFree} and generic irreducibility of the slowness polynomial~\cite[Theorem C]{dehoop2023reconstruction}.
Our proof is much simpler as squarefreeness is a weaker condition than irreducibility, and we believe the ideas may be useful to others pursuing similar results.
Let
\begin{equation}\label{eq:Msf}
M^{\mathrm{sf}}=\{m\in M_\C:\ \rho\tilde\rho\det\Gamma_c(e_n)\det\Gamma_{\tilde c}(e_n)\neq0,\ D(\rho,c,\,\cdot\,)\not\equiv0,\ D(\tilde\rho,\tilde c,\,\cdot\,)\not\equiv0\}.
\end{equation}
Its complement is defined by the vanishing of $\rho\tilde\rho\det\Gamma_c(e_n)\det\Gamma_{\tilde c}(e_n)$, or of all coefficients of $D(\rho,c,\,\cdot\,)$, or of all coefficients of $D(\tilde\rho,\tilde c,\,\cdot\,)$, so $M^{\mathrm{sf}}$ is Zariski open, and it contains $M^\mathrm{simp}$. If $m\in M^{\mathrm{sf}}$ and $D(\rho,c,\xi')\neq0$, all roots of $z\mapsto\det p(\xi',z)$ are simple, so $(\rho,c)$ is squarefree by Lemma~\ref{lem:simple-roots-sqf} with $j=n$, and likewise for $(\tilde\rho,\tilde c)$. Hence every $m\in M^{\mathrm{sf}}$ is squarefree. Moreover, $M^{\mathrm{sf}}\subset M^\circ$ with $M^\circ$ from the proof of Lemma~\ref{lem:generic-noncontainment}. Finally, $M^{\mathrm{sf}}$ is invariant under rotations of~$\R^n$ fixing~$e_n$. Let $\operatorname{Rot}_0\in O(n)$ be a rotation of $\R^n$ and denote by $\operatorname{Rot}_2 \in GL(\operatorname{Sym}^2(\operatorname{Sym}^2(\R^2)))$ the rotation of the stiffness tensor (see~\cite[Sec. 2]{Saksala26} for more details of how $\operatorname{Rot}_0$ and $\operatorname{Rot}_2$ are introduced in $n=2$). The action of $\operatorname{Rot}_2$ on a stiffness tensor $c$ are: $(\operatorname{Rot}_2 c)_{ijk\ell}=(\operatorname{Rot}_0)_{ii'}(\operatorname{Rot}_0)_{jj'}(\operatorname{Rot}_0)_{kk'}(\operatorname{Rot}_0)_{\ell\ell'}c_{i'j'k'\ell'}$. We write $\operatorname{Rot}_2 m=(\rho,\operatorname{Rot}_2 c,\tilde\rho,\operatorname{Rot}_2\tilde c)$ such that we have $\Gamma_{\operatorname{Rot}_2 c}(\xi)=\operatorname{Rot}_0\,\Gamma_c(\operatorname{Rot}_0^T\xi)\operatorname{Rot}_0^T$, hence $\det\Gamma_{\operatorname{Rot}_2 c}(e_n)=\det\Gamma_c(e_n)$ and $D(\rho,\operatorname{Rot}_2 c,\xi')=D(\rho,c,\operatorname{Rot}_0^T\xi')$, where $\operatorname{Rot}_0^T$ acts on $(\xi_1,\dots,\xi_{n-1})$.

\begin{proof}[Proof of Lemma~\ref{lem:alg-alessandrini}]
After a rotation that transforms the elastic operators and the DN maps covariantly and replaces $c$, $\tilde c$ by the rotated tensors, we may assume that the inner unit normal at~$x$ is~$e_n$, and we choose the boundary normal coordinates of Proposition~\ref{prop:Alessandrini} so that their differential at the boundary point is the identity. Then the covectors and polarizations in Proposition~\ref{prop:Alessandrini} at $x'=\hat x'$ are Euclidean ones, and $m=(\rho,\tilde\rho,c,\tilde c)$ denotes the values of the (rotated) coefficients at~$x$. Let $\hat x=(\hat x',0)\in(0,T)\times\partial\Omega$ correspond to the boundary point~$x$ in the statement of this Lemma. The values of $(\rho,c)$ and $(\tilde\rho,\tilde c)$ at~$\hat x$ are real, satisfy \eqref{eq:rhopos}--\eqref{eq:pos} and are squarefree by assumption. Let $p(\xi',z)$ and $\tilde p(\xi',z)$ denote the principal symbols for these values, and let $\delta>0$ be such that~\eqref{eq:pos} holds for both~$c$ and~$\tilde c$ at~$\hat x$. By~\eqref{eq:Gamma-elliptic} we have $\det\Gamma_c(e_n)>0$ and $\det\Gamma_{\tilde c}(e_n)>0$, and for real $\xi'=(\xi_0,\dots,\xi_{n-1})$ with $\delta(\xi_1^2+\dots+\xi_{n-1}^2)>2\max\{\rho,\tilde\rho\}\xi_0^2$ the matrices $p(\xi',z)$ and $\tilde p(\xi',z)$ are positive definite for all real~$z$, so $\det p(\xi',\,\cdot\,)$ and $\det\tilde p(\xi',\,\cdot\,)$ have no real roots for such $\xi'$. These $\xi'$ form a non-empty open set $\mathcal E_0\subset\R^n$. As $\det\Gamma_c(e_n)\neq0$, the polynomial $\det p(\xi',\,\cdot\,)$ has degree $2n$ for every~$\xi'$, and its discriminant $D(\rho,c,\xi')$ from the proof of Lemma~\ref{lem:generic-simplicity} vanishes exactly when it has a multiple root; by Lemma~\ref{lem:squareFree} and squarefreeness, the polynomial $D(\rho,c,\,\cdot\,)D(\tilde\rho,\tilde c,\,\cdot\,)$ is therefore not identically zero. Its zero set in~$\R^n$ has empty interior, since a polynomial vanishing on a non-empty open subset of~$\R^n$ vanishes identically, so there is $\hat\xi'\in\mathcal E_0$ at which it does not vanish. Then $\det p(\hat\xi',\,\cdot\,)$ and $\det\tilde p(\hat\xi',\,\cdot\,)$ each have $2n$ simple non-real roots; as their coefficients are real, the roots come in conjugate pairs, so exactly~$n$ roots $z^{(r)}$, respectively $\tilde z^{(s)}$, have positive imaginary part. In particular $(\hat\xi',z^{(r)})\in\Xi_1(P,\hat x)$ and $(\hat\xi',\tilde z^{(s)})\in\Xi_1(\tilde P,\hat x)$ for all~$r,s$, so the roots and kernel vectors depend smoothly on $(x',\xi')$ near $(\hat x',\hat\xi')$ by Lemma~\ref{lemma_z_omega}, and all hypotheses of Proposition~\ref{prop:Alessandrini} are satisfied.

Fix one root $z=z^{(r)}$ and one root $\tilde z=\tilde z^{(s)}$, with the kernel vectors $\omega$ and $\tilde\omega$ from Lemma~\ref{lemma_z_omega}, and restrict to $x'=\hat x'$, so that the coefficients are the values~$m$ at~$\hat x$. For real~$\xi'$ the polynomial $\det\tilde p(\xi',\,\cdot\,)$ has real coefficients, so $\overline{\tilde z(\xi')}$ is again a simple root of it, and $\tilde p(\xi',\overline{\tilde z})\,\overline{\tilde\omega}=\overline{\tilde p(\xi',\tilde z)\,\tilde\omega}=0$. Hence, for real~$\xi'$ near~$\hat\xi'$, the point
\begin{equation}\label{eq:bridge-point}
P(\xi')=\big((\xi',z(\xi')),\,(\xi',\overline{\tilde z(\xi')}),\,\omega(\xi'),\,\overline{\tilde\omega(\xi')}\big)
\end{equation}
lies in $A_\R\cap W_m$, and since $c$, $\tilde c$, $\rho$, $\tilde\rho$ and~$\xi_0$ are real, the identity of Proposition~\ref{prop:Alessandrini} at $x'=\hat x'$ reads $F(P(\xi'))=0$. Thus the algebraic covectors are $\xi=\zeta$ and $\tilde\xi=\overline{\tilde\zeta}$, the algebraic polarizations are $\omega$ and $\overline{\tilde\omega}$, and all complex conjugates are absorbed in the choice of the point~$P(\xi')$; the polynomials $F$, $f_i$ and~$\tilde f_i$ contain none.

These points fill an open subset of $A_\R\cap W_m$. Indeed, work in an affine chart of~$\CP^{n-1}$ around~$[\hat\xi']$ given by a real coordinate. By the holomorphic implicit function theorem, the simple roots $z(\hat\xi')$ and $\overline{\tilde z(\hat\xi')}$ extend to holomorphic functions of complex~$\xi'$ near~$\hat\xi'$, and near $(\hat\xi',z(\hat\xi'))$, respectively $(\hat\xi',\overline{\tilde z(\hat\xi')})$, every zero of $\det p$, respectively $\det\tilde p$, is of this form. At simple roots the kernels are one-dimensional and spanned by a non-zero column of the adjugate, so the projective kernel vectors extend holomorphically as well. Hence, near $P(\hat\xi')$, the variety~$W_m$ is the image of a holomorphic extension of $\xi'\mapsto P(\xi')$ to complex~$\xi'$, and $A_\R\cap W_m$ is the image of the real~$\xi'$. For a small real neighbourhood~$V$ of~$\hat\xi'$, the set $O=\{P(\xi'):\xi'\in V\}$ is therefore a non-empty open subset of $A_\R\cap W_m$ with $O\subset U_m$, at whose points $\xi_n$ and~$\tilde\xi_n$ are simple roots. Hence $m$ is a valid pair.
\end{proof}

\begin{lemma}[Smoothness]\label{lem:smooth}
Fix any $m\in M_\C$.
Let $x=(\xi,\tilde\xi,\omega,\tilde\omega)\in W_m$ be such that~$\xi_n$ is a simple root of $\xi_n\mapsto\det(p(\xi))$ and similarly with tildes.
Then~$x$ is a smooth point of~$W_m$.
\end{lemma}

\begin{proof}
Under the assumption of simple roots the kernels of the symbols are one-dimensional.

In the locus where $\dim(\ker(p))=1$, the vector $\omega\in\ker(p)$ depends smoothly and uniquely on~$p$.
Therefore we may simplify eliminate the variables~$\omega$ and~$\tilde\omega$ (which were not assumed real) and study the smooth points of the resulting variety.

Eliminating $\omega\in\CP^{n-1}$ from the~$n$ equations $p\omega=0$ leads (famously) to the one equation $\det(p)=0$.
Therefore we end up with the variety
\begin{equation}
Z
=
\{
(\xi,\tilde\xi)
\in
\Xi
;
\det(p(\xi))
=
\det(\tilde p(\tilde\xi))
=
0
\}
,
\end{equation}
whose smoothness we need to establish.

This~$Z$ is a smooth codimension~$2$ (i.e. dimension $n-1$) subvariety of~$\Xi$ wherever the gradients of the determinants are linearly independent.
Fortunately~$\det(p(\xi))$ is independent of $\tilde\xi_n$ and $\det(\tilde p(\tilde\xi))$ of $\tilde\xi_n$ and by assumption
\begin{equation}
\frac{\partial\det(p(\xi))}{\partial\xi_n}
\neq
0
\neq
\frac{\partial\det(\tilde p(\tilde\xi))}{\partial\tilde\xi_n}
.
\end{equation}
This establishes the required smoothness.
\end{proof}

\begin{lemma}[Real fibrations]\label{lem:fibrations}
Let $X,Y$ be complex varieties, where~$X$ is defined over~$\R$, and let $W\subset X\times Y$ be an irreducible
subvariety. Denote by $\pi:W\to X$ the projection.
Suppose that~$O$ is a non-empty Euclidean open subset
of $W\cap(X(\R)\times Y)$ containing a point
$p=(x,y)$ such that~$X$ is non-singular at~$x$ and~$\pi$ is smooth at~$p$. Then~$O$ is Zariski dense in~$W$.
\end{lemma}

\begin{proof}
Since~$X$ is non-singular at the real point~$x$, we may
choose local holomorphic coordinates
$(z_1,\ldots,z_r)$ on~$X$ near~$x$ such that~$X(\R)$ is locally given by
$\im(z_1)=\cdots=\im(z_r)=0$.
As~$\pi$ is smooth at~$p$, its analytification is a
holomorphic submersion near~$p$. By the holomorphic
implicit function theorem, there are local holomorphic
coordinates
\begin{equation}
(z_1,\ldots,z_r,w_1,\ldots,w_d)
\end{equation}
on~$W$ near~$p$ in which~$\pi$ is the projection onto
the first~$r$ coordinates. In particular, the subset
$W\cap(X(\R)\times Y)$ is locally identified with
$\R^r\times\C^d$.

Let~$f$ be a regular function on an affine open subset
of~$W$ containing~$p$ which vanishes on~$O$.
Since~$O$ is open in the Euclidean topology,
we may shrink the coordinate neighborhood of~$p$ so
that~$f$ vanishes on a set of the form
\begin{equation}
V_\R\times V_\C\subset\R^r\times\C^d,
\end{equation}
where~$V_\R$ and~$V_\C$ are non-empty open subsets.
For every fixed $w\in V_\C$, the holomorphic function
$z\mapsto f(z,w)$ vanishes on~$V_\R$.
Since the imaginary partial derivatives of~$f$ are uniquely defined by the real partial derivatives in the $z$-coordinates, Taylor expanding~$f$ on~$V_\R$ shows that~$f$ vanishes
on a non-empty complex open neighborhood of~$p$ in~$W$.

It follows that the zero set of~$f$ contains a non-empty
Euclidean open subset of~$W$. Since~$W$ is irreducible,
a proper algebraic subset of~$W$ cannot contain such
an open subset. Thus~$f$ vanishes identically on~$W$.
This proves that the Zariski closure of~$O$ is~$W$.
\end{proof}

The main result (that goes well with the analytic result) would be:

\begin{theorem}[Generic uniqueness]\label{thm:generic-uniqueness}
For a valid pair $m = ((\rho,c),(\tilde \rho, \tilde c))\in G\cap M^{\mathrm{sf}}$%

we have $\rho=\tilde\rho$ and $c=\tilde c$.
\end{theorem}

\begin{proof}
Let $m\in G\cap M^{\mathrm{sf}}$ be a valid pair. Since $M^{\mathrm{sf}}\subset M^\circ$, the variety $W_m$ lies in the open set~$B$ from the proof of Lemma~\ref{lem:generic-noncontainment}.

Let $O\subset A_\R\cap W_m\cap U_m$ be the non-empty
Euclidean open subset appearing in the definition of
a valid pair. Consider the projection
\begin{equation}
\pi_m:W_m\longrightarrow\CP^{n-1},
\qquad
(\xi,\tilde\xi,\omega,\tilde\omega)
\longmapsto[\xi_0:\cdots:\xi_{n-1}].
\end{equation}
On the open set $B$ this projection is well-defined,
and its restriction to $A_\R\cap W_m$ takes values
in~$\RP^{n-1}$.

Choose an affine chart of $\CP^{n-1}$ containing
$[\xi_0:\cdots:\xi_{n-1}]$ and write $x$ for its
coordinates, $y=\xi_n$ and $z=\tilde\xi_n$,
after normalizing the common homogeneous coordinates.
The equations defining~$W_m$ take the form
\begin{equation}
p(x,y)\omega=0,\qquad
\tilde p(x,z)\tilde\omega=0.
\end{equation}

By the definition of a valid pair, $O$ contains a point
$p=(\xi,\tilde\xi,\omega,\tilde\omega)$ at which
$\xi_n$ and $\tilde\xi_n$ are simple roots of
the corresponding determinant equations $\det p(x,y) = 0$ and $\det \tilde p(x,z) = 0$.
We show that $\pi_m$ is smooth at $p$.
By the choice of~$p$,
\begin{equation}
\partial_y\det p(x,y)\ne0,\qquad
\partial_z\det\tilde p(x,z)\ne0
\end{equation}
at~$p$. In particular, both matrices have rank $n-1$.
Indeed, if either matrix had rank at most $n-2$,
its adjugate would vanish and Jacobi's formula for derivatives of determinants of matrix valued functions would
give that the corresponding derivative of its
determinant is zero.

By the holomorphic implicit function theorem, the
two determinant equations determine $y$ and $z$
as holomorphic functions of $x$ near $p$. Since the
kernels are constantly one-dimensional, the projective kernel
vectors $\omega$ and $\tilde\omega$ also depend
holomorphically on $x$. Thus $\pi_m$ is locally
a biholomorphic map near $p$, and in particular is a
smooth morphism at~$p$.

Hence $W_m$ is smooth at~$p$, so $p$ lies on a unique irreducible component~$C$ of~$W_m$, and $W_m$ coincides with~$C$ near~$p$. Then $O\cap C$ contains a non-empty
Euclidean open subset of
\begin{equation}
C\cap\pi_m^{-1}(\RP^{n-1}).
\end{equation}
By Lemma~\ref{lem:fibrations}, this subset is Zariski
dense in $C$. Since $O\subset U_m$ and $U_m$ is
Zariski closed, it follows that
\begin{equation}
C\subset U_m.
\end{equation}
This contradicts the assumption $m\in G$.
\end{proof}

\begin{proof}[Proof of Theorem~\ref{thm:main}]
For a rotation~$\operatorname{Rot}_0$ of~$\R^n$ fixing~$e_n$, the map
\begin{equation}
(\xi,\tilde\xi,\omega,\tilde\omega)\mapsto(\operatorname{Rot}_0\xi,\operatorname{Rot}_0\tilde\xi,\operatorname{Rot}_0\omega,\operatorname{Rot}_0\tilde\omega)
,
\end{equation}
with $\operatorname{Rot}_0$ acting on the spatial components, preserves~$\Xi$ and~$A_\R$ and maps $W_m$ onto $W_{\operatorname{Rot}_2 m}$ and $U_m$ onto $U_{\operatorname{Rot}_2 m}$, since $F$, $f_i$ and~$\tilde f_i$ are built from contractions of tensors. Hence~$G$ is invariant under such rotations, and so is~$M^{\mathrm{sf}}$. Let $K_0\subset G$ be the non-empty Zariski open subset from Lemma~\ref{lem:generic-noncontainment}, and let~$K$ be the union of the sets $\operatorname{Rot}_{2}(K_0\cap M^{\mathrm{sf}})\cap M_\R$ over all rotations~$\operatorname{Rot}_{0}$ fixing~$e_n$. Then $K\subset G\cap M^{\mathrm{sf}}$, $K$ is invariant and open in~$M_\R$, and it is dense in~$M_\R$: it contains the real points of $K_0\cap M^{\mathrm{sf}}$, whose complement is a proper algebraic subset of the irreducible space~$M_\C$, and a non-zero polynomial cannot vanish on a non-empty open subset of~$M_\R$.

Let $\operatorname{Rot}_{0,p}$ be a rotation mapping the inner unit normal at~$p$ to~$e_n$, and set
\begin{equation}
G_p=\{m\in(\R_+\times\possymsym{\R^n})^2:\ \operatorname{Rot}_{2,p} m\in K\}.
\end{equation}
This set is open and dense in $(\R_+\times\possymsym{\R^n})^2$, and by the invariance of~$K$ it does not depend on the choice of~$\operatorname{Rot}_{0,p}$.

Let $m(p)=(\rho(p),c(p),\tilde\rho(p),\tilde c(p))\in G_p$ and $\Lambda_{\rho,c}=\Lambda_{\tilde\rho,\tilde c}$. After the rotation by~$\operatorname{Rot}_{0,p}$, the inner unit normal at~$p$ is~$e_n$, and the values of the coefficients at~$p$ are $\operatorname{Rot}_{2,p} m(p)\in K$, which is squarefree. By Lemma~\ref{lem:alg-alessandrini}, $\operatorname{Rot}_{2,p} m(p)$ is a valid pair, and by Theorem~\ref{thm:generic-uniqueness}, $\rho(p)=\tilde\rho(p)$ and $\operatorname{Rot}_{2,p} c(p)=\operatorname{Rot}_{2,p} \tilde c(p)$, that is, $c(p)=\tilde c(p)$.
\end{proof}

\section{Stiffness tensor fields}
\label{sec:fields}

Our last result promotes the pointwise result to one for stiffness tensor fields.
It is based on a simple topological fact.

\begin{lemma}
\label{lma:section-to-field}
Let~$M$ be a smooth closed manifold and~$E$ a real vector bundle atop it, both of finite dimension.
Suppose there is a closed subset $S\subset E^2$ with nowhere dense fibers.
Then for a residual (i.e., comeagre) set of pairs of sections $(f_1,f_2)\in C^\infty(M;E)^2$ the set $\{x\in M;(f_1(x),f_2(x))\in S_x\}$ has no interior points.
\end{lemma}

\begin{proof}
Take a countable basis $(U_i)_{i\in\N}$ for the topology of~$M$ with local trivializations of~$E$.
For each~$i$ denote by $A_i\subset C(M;E)^2$ the set of pairs $(f_1,f_2)$ for which $(f_1(x),f_2(x))\notin S_x$ for some $x\in U_i$.

Each~$A_i$ is open because if $(f_1(x),f_2(x))\notin S_x$, then by closedness of~$S_x$ this is also true at~$x$ for a neighborhood of $(f_1,f_2)$.

Each~$A_i$ is dense because if $(f_1,f_2)\notin A_i$, then adding a small bump function supported in~$U_i$ makes $(f_1(x),f_2(x))\notin S_x$ for some $x\in U_i$.
Any point~$x$ where the bump function is non-zero works because the fiber~$S_x$ above it has no interior points.

The space $C^\infty(M;E)$ is Fréchet.
Therefore by the Baire category theorem the desired conclusion holds in the residual set $\cap_{i\in\N}A_i$.
\end{proof}

Our fibers are much smaller than ``closed without interior points'', but in general a similar statement for pairs of sections does not hold in an open set.
The statement on the level of fields is therefore naturally weaker.

\begin{proof}[Proof of Corollary~\ref{cor:fields}]
By continuity of the density and stiffness tensor fields it is enough to prove equality in a dense subset of~$\partial\Omega$.
We identify~$\R$ with~$\R_+$ through the exponential map, and similarly for positive definite stiffness tensors.
With this identification we get to study all smooth (pairs of) sections of a vector bundle.
We may thus apply Lemma~\ref{lma:section-to-field} for $M=\partial\Omega$ and $E=\R_+\times\possymsym{\R^n}\approx\R^N$ for some $N\in\N$.
The fiber~$S_x$ is the complement of the set~$G_x$ from the proof of Theorem~\ref{thm:main}, which is closed with empty interior. Since~$K$ is invariant under rotations fixing~$e_n$, $S_x$ depends on $x\in\partial\Omega$ only through the normal direction at~$x$, and choosing the rotations~$\operatorname{Rot}_{0,x}$ continuously on small open subsets of~$\partial\Omega$ shows that $S\subset E^2$ is closed in the bundle.
All assumptions are satisfied, so the claim follows.
\end{proof}

\printbibliography

@article{Nakamura1993,
  title = {Identification of Lame Parameters by Boundary Measurements},
  volume = {115},
  ISSN = {0002-9327},
  url = {http://dx.doi.org/10.2307/2375069},
  DOI = {10.2307/2375069},
  number = {5},
  journal = {American Journal of Mathematics},
  publisher = {JSTOR},
  author = {Nakamura,  Gen and Uhlmann,  Gunther},
  year = {1993},
  month = Oct,
  pages = {1161}
}

@incollection{SU90,
  author    = {Sylvester, John and Uhlmann, Gunther},
  title     = {Inverse problems in anisotropic media},
  booktitle = {Inverse Scattering and Applications (Amherst, MA, 1990)},
  series    = {Contemp. Math.},
  volume    = {122},
  pages     = {105--117},
  publisher = {Amer. Math. Soc.},
  address   = {Providence, RI},
  year      = {1991}
}

@article{SU05,
  author  = {Stefanov, Plamen and Uhlmann, Gunther},
  title   = {Boundary rigidity and stability for generic simple metrics},
  journal = {J. Amer. Math. Soc.},
  volume  = {18},
  pages   = {975--1003},
  year    = {2005}
}

@article{Mo14,
  author  = {Montalto, Carlos},
  title   = {Stable determination of a simple metric, a covector field and a potential from the hyperbolic {D}irichlet-to-{N}eumann map},
  journal = {Comm. Partial Differential Equations},
  volume  = {39},
  pages   = {120--145},
  year    = {2014}
}

@article{SY18,
  author  = {Stefanov, Plamen and Yang, Yang},
  title   = {The inverse problem for the {D}irichlet-to-{N}eumann map on {L}orentzian manifolds},
  journal = {Anal. PDE},
  volume  = {11},
  pages   = {1381--1414},
  year    = {2018}
}

@misc{Saksala26,
  doi = {10.48550/ARXIV.2610.13747},
  author = {Ilmavirta,  Joonas and Saksala,  Teemu and Schl\"{u}ter,  Hjørdis},
  title = {Travel time tomography for elliptically anisotropic elastic planar media},
  publisher = {arXiv},
  year = {2026}
}

@article{KV84,
  author  = {Kohn, Robert and Vogelius, Michael},
  title   = {Determining conductivity by boundary measurements},
  journal = {Comm. Pure Appl. Math.},
  volume  = {37},
  pages   = {289--298},
  year    = {1984}
}

@article{SU88,
  author  = {Sylvester, John and Uhlmann, Gunther},
  title   = {Inverse boundary value problems at the boundary---continuous dependence},
  journal = {Comm. Pure Appl. Math.},
  volume  = {41},
  pages   = {197--219},
  year    = {1988}
}

@article{LU89,
  author  = {Lee, John M. and Uhlmann, Gunther},
  title   = {Determining anisotropic real-analytic conductivities by boundary measurements},
  journal = {Comm. Pure Appl. Math.},
  volume  = {42},
  pages   = {1097--1112},
  year    = {1989}
}

@article{DKSU09,
  author  = {Dos Santos Ferreira, David and Kenig, Carlos E. and Salo, Mikko and Uhlmann, Gunther},
  title   = {Limiting {C}arleman weights and anisotropic inverse problems},
  journal = {Invent. Math.},
  volume  = {178},
  pages   = {119--171},
  year    = {2009}
}

@phdthesis{stolk00,
    author =   {Stolk,Christiaan Cornelis},
    title = {On the modeling and inversion of seismic data},
    school = {Utrecht University},
    year = {2000}
}

@book{Smith2000Invitation,
  title     = {An Invitation to Algebraic Geometry},
  author    = {Smith, Karen~E. and Kahanp{\"a}{\"a}, Lauri and Kek{\"a}l{\"a}inen, Pekka and Traves, William},
  series    = {Universitext},
  publisher = {Springer New York},
  year      = {2000},
  doi       = {10.1007/978-1-4757-4497-2}
}

@book{atiyah1969introduction,
  title     = {Introduction to Commutative Algebra},
  author    = {Atiyah, M. F. and MacDonald, I. G.},
  year      = {1969},
  publisher = {Addison-Wesley},
  address   = {Reading, Massachusetts},
}

@book{Hartshorne,
  title={Algebraic {G}eometry},
  author={Hartshorne, R.},
  isbn={9780387902449},
  lccn={lc77001177},
  series={Graduate Texts in Mathematics},
  url={https://books.google.at/books?id=3rtX9t-nnvwC},
  year={1977},
  publisher={Springer}
}

@misc{stacks-project,
  author       = "{The Stacks Project Authors}",
  title        = "Stacks Project",
  howpublished = "\url{https://stacks.math.columbia.edu/}"
}

@article{Oksanen2024,
  title = {Inverse problems for real principal type operators},
  volume = {146},
  ISSN = {1080-6377},
  url = {http://dx.doi.org/10.1353/ajm.2024.a917541},
  DOI = {10.1353/ajm.2024.a917541},
  number = {1},
  journal = {American Journal of Mathematics},
  publisher = {Project MUSE},
  author = {Oksanen, Lauri and Salo, Mikko and Stefanov, Plamen and Uhlmann, Gunther},
  year = {2024},
  month = Feb,
  pages = {161--240}
}

@article{Rachele001,
author = { Lizabeth   V Rachele },
title = {Boundary determination for an inverse problem in elastodynamics},
journal = {Communications in Partial Differential Equations},
volume = {25},
number = {11-12},
pages = {1951-1996},
year  = {2000},
publisher = {Taylor & Francis},
doi = {10.1080/03605300008821575},
}

@article{HoopNakamuraZhai19,
author = {de Hoop, Maarten V. and Nakamura, Gen and Zhai, Jian},
title = {Unique Recovery of Piecewise Analytic Density and Stiffness Tensor from the Elastic-Wave Dirichlet-To-Neumann Map},
journal = {SIAM Journal on Applied Mathematics},
volume = {79},
number = {6},
pages = {2359-2384},
year = {2019},
doi = {10.1137/18M1232802},
}

@article{Ilmavirta2025,
  title = {Gauge freedoms in the anisotropic elastic Dirichlet-to-Neumann map},
  volume = {19},
  ISSN = {1930-8345},
  url = {http://dx.doi.org/10.3934/ipi.2024055},
  DOI = {10.3934/ipi.2024055},
  number = {4},
  journal = {Inverse Problems and Imaging},
  publisher = {American Institute of Mathematical Sciences (AIMS)},
  author = {Ilmavirta,  Joonas and Schl\"{u}ter,  Hjørdis},
  year = {2025},
  pages = {816–828}
}

@misc{dehoop2023reconstruction,
      title={Reconstruction of generic anisotropic stiffness tensors from partial data around one polarization}, 
      author={Maarten V. de Hoop and Joonas Ilmavirta and Matti Lassas and Anthony Várilly-Alvarado},
      year={2023},
      eprint={2307.03312},
      archivePrefix={arXiv},
      primaryClass={math.DG}
}

@unpublished{maxwell-ag,
    author = {Antonio Cocan and Maarten V. de Hoop and Joonas Ilmavirta and Matti Lassas and Anthony V\'{a}rilly-Alvarado},
    title = {{Generic Recovery of Permittivity and Permeability in Anisotropic Maxwell Systems}},
    month = aug,
    year = 2026,
    note = {arXiv:2608.02900}
}

@unpublished{elastic-regularity,
    author = {Joonas Ilmavirta and Pieti Kirkkopelto and Antti Kykk\"{a}nen},
    title = {{Horizontal and Vertical Regularity of Elastic Wave Geometry}},
    month = nov,
    year = 2025,
    note = {arXiv:2511.16466}
}

@phdthesis{stolk-thesis,
    author = {C. C. Stolk},
    title = {On the modeling and inversion of seismic data},
    school = {Utrecht University},
    year = {2000}
}

\end{document}